\documentclass[12pt]{elsarticle}

\usepackage{amsmath}
\usepackage{amsfonts} 
\usepackage{amssymb}  
\usepackage{graphicx}
\usepackage{amsthm} 
\usepackage{enumerate}
\usepackage{caption}
\usepackage{subfigure}

\allowdisplaybreaks[4] 

\numberwithin{equation}{section}
\newtheorem{thm}{Theorem}[section] 
\newtheorem{lem}{Lemma}[section]
\newtheorem{rmk}{Remark}[section]
\newtheorem{prop}{Proposition}[section]

\vspace{0.4cm}

\begin{document}

\title
{\textbf{
A proof of the uniqueness of the limit cycle of a
quasi-homogeneous system}}

\author[1]{Ziwei Zhuang}
\ead{zhuangzw@mail2.sysu.edu.cn}
\author[1]{Changjian Liu}
\ead{liuchangj@mail.sysu.edu.cn}

\address[1]{School of Mathematics (Zhuhai),
Sun Yat-sen University, 519086,
Zhuhai, P. R. China}

\begin{abstract}
A. Gasull shared a list of 33
open problems in low dimensional dynamical systems
in his work in 2021.
The second part of Problem 3 is about whether the limit
cycle of a quasi-homogeneous system
$
  \dot{x}=y,\;
  \dot{y}=-x^3+\alpha x^2y+y^3
$
is unique.
In this paper,
we give a positive answer to
this question
by analysing the uniqueness of
the heteroclinic separatrix at infinity.

\end{abstract}


\maketitle

\section{Introduction}
Consider the uniqueness of the limit cycle
 of the following concrete system with a real parameter
 $\alpha:$
\begin{equation}\label{concretesys}
\left\{
       \begin{aligned}
       &\frac{\mathrm{d}x}{\mathrm{d}t}=y,\\
       &\frac{\mathrm{d}y}{\mathrm{d}t}=-x^3+\alpha x^2y+y^3.
       \end{aligned}
\right.
\end{equation}
This problem is put forward by A. Gasull in
\cite[Problem 3(ii)]{gasull2021some}.
The origin of the system, which is the unique and
fixed equilibrium,
is a stable focus for $\alpha<0,$
while it is an unstable focus for $\alpha\ge0,$
see \cite[Corollary 5]{cima1995cyclicity}.
Since for any $xy\neq0,$
$$
\left|
\begin{array}{ll}
y & -x^3+\alpha x^2y+y^3 \\
\frac{\partial}{\partial \alpha}y
& \frac{\partial}{\partial \alpha} (-x^3+\alpha x^2y+y^3)
\end{array}
\right|
=x^2y^2>0
$$
and
$$
\frac{\mathrm{d}y}{\mathrm{d}x}
=\frac{y^3-x^3}{y}+\alpha x^2
\rightarrow\pm\infty\quad
\text{as}\quad
\alpha\rightarrow\pm\infty,
$$
the system is a
semi-complete family of rotated vector fields
(mod $xy=0$) with respect to $\alpha,$ which is defined in
\cite[Definitions 2, 3]{perko1993rotated}.
From the properties of rotated vector fields,
an unstable limit cycle bifurcates from the origin for
$\alpha<0$ and $\alpha$ sufficiently close to 0,
which is also a generalized Hopf bifurcation.
For the number of the limit cycle of the system,
the following results are already known.

\begin{lem}[\cite{gasull2006upper}]
\label{lemma_Gasull2006}
(i)
System \eqref{concretesys} has no limit cycles when
$\alpha\ge0$ or $\alpha<-2.679.$

(ii)
System \eqref{concretesys} has at most one limit
cycle when $0>\alpha>-3/\sqrt[3]{2}\approx-2.381.$
The limit cycle is hyperbolic and unstable if
it exists.
\end{lem}

\begin{lem}[\cite{giacomini2015transversal}]
\label{lemma_Giacomini2015}
System \eqref{concretesys} has at least one
limit cycles when
$0>\alpha>-3\sqrt[3]{6\sqrt3-9}/\sqrt[3]4\approx-2.1103.$
\end{lem}
Besides, it is also mentioned in
\cite{gasull2006upper} that
the existence range of limit cycles seems to be
$(-2.198,0)$ by a numerical computation.

In this paper, we first analyse the behavior
of system \eqref{concretesys} at infinity by
the Poincar\'{e} transformation.
As one of the most important results,
we estimate the behavior of the separatrices
of two saddles at
$\alpha=-3/\sqrt[3]{2}\approx-2.3811.$
Then with the help of known results,
we show that there is a unique
$\alpha^*\in(-2.3811,-2.1103]$ such that
system \eqref{concretesys} has a
heteroclinic separatrix at infinity.
Furthermore, the system
has exactly one limit cycle for any
$\alpha\in(\alpha^*,0)$
and no cycles for the else region.

\section{Phase portrait at infinity}
We first study the phase portrait of system \eqref{concretesys} at infinity.
It is not hard to see that
$(\pm\infty,0)$
are not critical points at infinity.
Then, by the Poincar\'{e} transformation
$$ x=\frac{v}{z},\quad y=\frac{1}{z}, $$
system \eqref{concretesys} is transformed to
\begin{equation}\label{transformedSys(v,z)}
\left\{
        \begin{aligned}
        &\frac{\mathrm{d}v}{\mathrm{d}\tau}=vf(v,\alpha)+z^2,\\
        &\frac{\mathrm{d}z}{\mathrm{d}\tau}=zf(v,\alpha),
        \end{aligned}
\right.
\end{equation}
where $\mathrm{d}\tau=\mathrm{d}t/z^2$ and
\begin{equation*}
  f(v,\alpha)=v^3-\alpha v^2-1.
\end{equation*}
It's easy to see that for any $\alpha<0$ and
$\mu\in\left(-1,f\left(2\alpha/3,\alpha\right)\right),$
the equation
$f(v,\alpha)=\mu$
has three different real roots for $v.$
In ascending order, they are denoted by
\begin{equation}\label{def_of_vmui}
  v^{\mu}_i(\alpha),\; i=1,2,3.
\end{equation}
$v^{\mu}_i(\alpha)$ is a continuous function
with respect to $\alpha$ in some suitable
range, and we have
$v^{\mu}_1(\alpha)<2\alpha/3<v^{\mu}_2(\alpha)
<0<v^{\mu}_3(\alpha).$
Since
$f\left(2\alpha/3,\alpha\right)
=-4\alpha^3/27-1$
monotonically decreases as
$\alpha$ decreases to $-\infty.$
it is not hard to verify the following result.

\begin{lem}\label{lemma_features_of_vmui}
As $\alpha$ decreases to $-\infty,$ $v^{\mu}_1(\alpha)$ monotonically decreases to
$-\infty,$ $v^{\mu}_2(\alpha)$ monotonically increases to $0$
and $v^{\mu}_3(\alpha)$ monotonically decreases to $0.$
\end{lem}

Before the last lemma in this section,
i.e., Lemma \ref{lemma_L-_above_L+},
we only consider the case $\mu=0,$
for which
$v_i^\mu(\alpha)$ is the abscissa of
the equilibrium of system \eqref{transformedSys(v,z)}.
Precisely,  when
$\alpha<-3/\sqrt[3]{4},$
system \eqref{transformedSys(v,z)}
has four equilibria, which are respectively denoted by
\begin{equation*}
  P_0=(0,0),\;
  P_1(\alpha)=\left(v_1^0(\alpha),0\right),\;
  P_2(\alpha)=\left(v_2^0(\alpha),0\right),\;
  P_3(\alpha)=\left(v_3^0(\alpha),0\right).
\end{equation*}
Throughout this section, we always assume
$\alpha<-3/\sqrt[3]{4}\approx-1.8899.$

\subsection{Behavior of the trajectories for $\alpha<-3/\sqrt[3]{4}$}
\begin{figure}[htbp]
\centering
\includegraphics[width=1\linewidth]{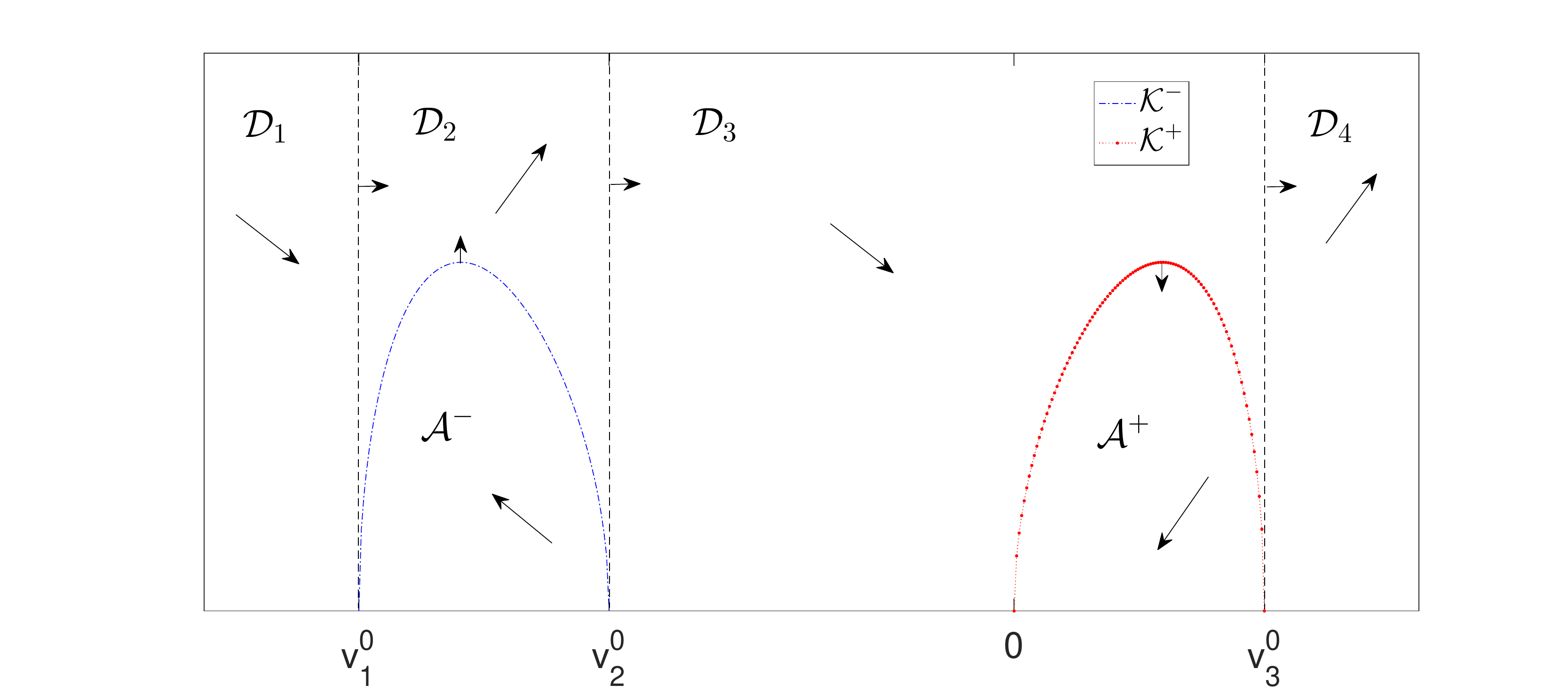}
\caption{Division of domains in the upper half plane
and approximate directions of tangent vectors in each domain.}
\label{fig:domains}
\end{figure}

Since system \eqref{proposition_equilibria_of_sys(v,z)} is symmetric
with respect to $v$-axis,
we only need to consider the upper half plane.
Divide the upper half plane into six domains
in each of which
$\mathrm{d}v/\mathrm{d}\tau$ and
$\mathrm{d}z/\mathrm{d}\tau$
do not change signs
(see Figure \ref{fig:domains}):
\begin{equation*}
\begin{aligned}
&\mathcal{D}_1(\alpha)=\left(-\infty,v_1^0(\alpha)\right]
\times\mathbb{R}^+,
\\
&\mathcal{D}_2(\alpha)=\left\{
(v,z)\in\mathbb{R}\times\mathbb{R}^+\;|
\;v_1^0(\alpha)<v<v_2^0(\alpha),
\;vf(v,\alpha)+z^2>0
\right\},
\\
&\mathcal{D}_3(\alpha)=\left\{
(v,z)\in\mathbb{R}\times\mathbb{R}^+
\;|\;v_2^0(\alpha)\le v< v_3^0(\alpha),
\;vf(v,\alpha)+z^2>0
\right\},
\\
&\mathcal{D}_4(\alpha)=\left[v_3^0(\alpha),+\infty\right)
\times\mathbb{R}^+,
\\
&\mathcal{A}^{-}(\alpha)=\left\{
(v,z)\in\mathbb{R}\times\mathbb{R}^+
\;|\;v_1^0(\alpha)<v<v_2^0(\alpha),
\;vf(v,\alpha)+z^2<0
\right\},
\\
&\mathcal{A}^{+}(\alpha)=\left\{
(v,z)\in\mathbb{R}\times\mathbb{R}^+
\;|\;0<v<v_3^0(\alpha),
\;vf(v,\alpha)+z^2<0
\right\},
\end{aligned}
\end{equation*}
and let
\begin{equation*}
  \mathcal{D}(\alpha)=\mathcal{D}_1(\alpha)
  \cup\mathcal{D}_2(\alpha)
  \cup\mathcal{D}_3(\alpha)
  \cup\mathcal{D}_4(\alpha),
  \quad
  \mathcal{A}(\alpha)=\mathcal{A}^{-}(\alpha)
  \cup\mathcal{A}^{+}(\alpha).
\end{equation*}
Denote the curves by
\begin{equation*}
\begin{aligned}
  &\mathcal{K}^{-}(\alpha)=\left\{
  \left(
  v,\sqrt{-vf(v,\alpha)}
  \right) \;|\;
  v_1^0(\alpha)<v<v_2^0(\alpha)
  \right\},
  \\
  &\mathcal{K}^{+}(\alpha)= \left\{
  \left(
  v,\sqrt{-vf(v,\alpha)}
  \right) \;|\;
  0<v<v_3^0(\alpha)
  \right\}.
\end{aligned}
\end{equation*}
We see that
$\mathrm{d}v/\mathrm{d}\tau>0$
in $\mathcal{D},$
$\mathrm{d}v/\mathrm{d}\tau<0$
in $\mathcal{A}$
and
$\mathrm{d}v/\mathrm{d}\tau=0$
on $\mathcal{K}^-$ and
$\mathcal{K}^+.$
Moreover,
$\mathcal{K}^{-}$ and $\mathcal{K}^{+}$
are two traversals of
system \eqref{transformedSys(v,z)}.
The tangent vectors on
$\mathcal{K}^{-}$
direct from
$\mathcal{A}^{-}$
to $\mathcal{D}_2,$
and the ones on
$\mathcal{K}^{+}$
direct from
$\mathcal{D}_3$ to
$\mathcal{A}^{+}$
(also see Figure \ref{fig:domains}).
Note that
$\mathcal{A}^-$ is bounded by
the traversal $\mathcal{K}^-,$
the equilibria $P_1$ and
$P_2$
and the trajectory $\overline{P_2P_1}.$
If a trajectory
intersects $\mathcal{K}^-$ at a point $P,$
then the negative part from $P$ of this trajectory
entirely lies in $\mathcal{A}^-.$
Similarly, if a trajectory
intersects $\mathcal{K}^+$ at a point $P,$
then the positive part from $P$ of this trajectory
entirely lies in $\mathcal{A}^+.$

\begin{prop}\label{proposition_equilibria_of_sys(v,z)}
Assume that $\alpha<-3/\sqrt[3]{4}.$
\begin{enumerate}[(i)]
\item
$P_0(\alpha)$ is a hyperbolic stable node.

\item
$P_1(\alpha)$ is a semi-hyperbolic saddle,
whose stable separatrices entirely lies
on $v$-axis and unstable separatrices are
tangent to $v=v_1^0(\alpha)$ at $P_1(\alpha).$
Moreover, the unstable separatrix in the
upper half plane has no common with
$\mathcal{D}_1(\alpha),$
$\mathcal{A}^-(\alpha)$
or
$\mathcal{K}^-(\alpha).$

\item
$P_2(\alpha)$ is a semi-hyperbolic stable node.

\item
$P_3(\alpha)$ is a semi-hyperbolic saddle,
whose unstable separatrices entirely lies
on $v$-axis and stable separatrices are
tangent to $v=v_3^0(\alpha)$ at $P_3(\alpha).$
Moreover, the stable separatrix in the
upper half plane has no common with
$\mathcal{D}_4(\alpha),$
$\mathcal{A}^+(\alpha)$
or
$\mathcal{K}^+(\alpha).$
\end{enumerate}
\begin{proof}
One can easily see that
$P_0$ is a hyperbolic equilibrium
and that the linearized matrix at
$P_0$ has
two negative eigenvalues.
From the Perron's Theorem,
$P_0$ is a stable node.

By a translation
$\nu=v-v_i^0$ with $i\ge1,$
$P_i$ is translated to the origin
and system \eqref{transformedSys(v,z)}
is transferred to
\begin{equation}\label{transformedSys(nu,z)}
\begin{aligned}
  &\frac{\mathrm{d}\nu}{\mathrm{d}\tau}
  =\left(3v_i^0+2\alpha\right) \left(v_i^0\right)^2 \nu
  +Q(\nu,z),
  \\
  &\frac{\mathrm{d}z}{\mathrm{d}\tau}
  =R(\nu,z),
\end{aligned}
\end{equation}
where
\begin{equation*}
  \begin{aligned}
  &Q(\nu,z)=\left[
  3\left(2v_i^0+\alpha\right)v_i^0
  \right]\nu^2
  +\left( 4v_i^0+\alpha \right) \nu^3
  +\nu^4
  +z^2,
  \\
  &R(\nu,z)=\left[
  \left(3v_i^0+2\alpha\right)v_i^0 \nu
  +\left(3v_i^0+\alpha\right) \nu^2
  +\nu^3
  \right] z.
  \end{aligned}
\end{equation*}
Since
$\left(3v_i^0+2\alpha\right) \left(v_i^0\right)^2 \neq 0,$
the origin of system \eqref{transformedSys(nu,z)}
is a semi-hyperbolic equilibrium introduced
in many books, see \cite{dumortier2006qualitative} for example.
From the Implicit Function Theorem,
the equation
$\left(3v_i^0+2\alpha\right) \left(v_i^0\right)^2 \nu
  +Q(\nu,z)=0$
has a unique and analytic solution
$$
\nu=\nu(z)
=-\frac{2z^2}{\left(3v_i^0+2\alpha\right)(v_i^0)^2}+O(z^3)
$$
in a small neighborhood of the origin.
Substituting it into $R(\nu,z),$ we have
\begin{equation*}
  R\left(\nu(z),z\right)=-\frac{2z^3}{v_i^0} + O(z^{5}).
\end{equation*}
From \cite[Theorem 2.19 and Remark 2.20]{dumortier2006qualitative},
there is an invariant analytic curve tangent to
the $\nu$-axis at the origin,
which is $\nu$-axis itself in our system.
Besides, the qualitative properties of the origin
is determined by the signs of the terms
$\left(3v_1^0+2\alpha\right) \left(v_1^0\right)^2$
and
$-2/v_1^0,$ which states as follows:
\begin{itemize}
  \item
  for $v_1^0,$ we have
  $\left(v_1^0+2\alpha\right) \left(v_1^0\right)^2<0$
  and
  $-2/v_1^0>0.$
  Then, the origin of system \eqref{transformedSys(nu,z)}
  is a topological saddle, whose unstable
  separatrices are tangent to the $z$-axis at the origin;

  \item
  for $v_2^0,$ we have
  $\left(v_2^0+2\alpha\right) \left(v_2^0\right)^2>0$
  and
  $-2/v_2^0>0.$
  Then, the origin of system \eqref{transformedSys(nu,z)}
  is an unstable topological node, whose trajectories
  near the origin except for
  the two lying on the $\nu$-axis are
  tangent to the $z$-axis at the origin;

  \item
  for $v_3^0,$ we have
  $\left(v_3^0+2\alpha\right) \left(v_3^0\right)^2>0$
  and
  $-2/v_3^0<0.$
  Then, the origin of system \eqref{transformedSys(nu,z)}
  is a topological saddle, whose stable
  separatrices are tangent to the $z$-axis at the origin.
\end{itemize}

Moreover, if any trajectory $L$ passes through a point
$(v_0,z_0)$
in $\mathcal{D}_1$ in which
$\mathrm{d}v/\mathrm{d}\tau>0$
and
$\mathrm{d}z/\mathrm{d}\tau\le0,$
the $\alpha$-limit set of $L$ is in
$\left\{v<v_0,z\ge z_0\right\}\subset \mathcal{D}_1$
which doesn't contain $P_1.$
On the other hand, recall that
$\mathcal{A}^-$ is bounded by
the traversal $\mathcal{K}^-,$
the equilibria $P_1$ and $P_2$
and the trajectory $\overline{P_2P_1}.$
If any trajectory $L$ passes through
a point on $\mathcal{K}^-$ or
in $\mathcal{A}^-$ in which
$\mathrm{d}v/\mathrm{d}\tau<0$
$\mathrm{d}z/\mathrm{d}\tau>0,$
$L$ must be negatively approaches $P_2.$
Therefore, the unstable separatrix of
$P_1$ in the upper half plane has no common
with $\mathcal{D}_1,\mathcal{A}^-$
or $\mathcal{K}^-.$
By analogous analysis,
the stable separatrix of
$P_3$ in the upper half plane has no common
with $\mathcal{D}_4,\mathcal{A}^+$
or $\mathcal{K}^+.$

\end{proof}
\end{prop}

Denote the solution of system \eqref{transformedSys(v,z)} by
\begin{equation*}
  L(v_0,z_0,\alpha):\;
  \left(v(\tau;v_0,z_0,\alpha),
  z(\tau;v_0,z_0,\alpha)\right),
\end{equation*}
or
$\left(v(\tau),z(\tau)\right)$
for short,
with initial condition
$v(0;v_0,z_0,\alpha)=v_0$
and
$z(0;v_0,z_0,\alpha)=z_0.$
Assume that $(v_0,z_0)\in\mathcal{D}(\alpha).$
Let
\begin{equation*}
\begin{aligned}
  \underline{\tau}(v_0,z_0,\alpha)
  =\inf\left\{\right.
  &\tau<0\,|\,
  v(\tau_1;v_0,z_0,\alpha)f
  \left(v(\tau_1;v_0,z_0,\alpha),\alpha\right)
  +\left(z(\tau_1;v_0,z_0,\alpha)\right)^2>0,
  \\
  &\left.\forall \tau_1\in[\tau,0]\right\},
  \\
  \overline{\tau}(v_0,z_0,\alpha)
  =\sup\left\{\right.
  &\tau>0\,|\,
  v(\tau_1;v_0,z_0,\alpha)f
  \left(v(\tau_1;v_0,z_0,\alpha),\alpha\right)
  +\left(z(\tau_1;v_0,z_0,\alpha)\right)^2>0,
  \\
  &\left.\forall \tau_1\in[0,\tau]\right\}.
\end{aligned}
\end{equation*}
Since
$\mathrm{d}v(\tau;v_0,z_0,\alpha)/\mathrm{d}\tau>0$
on $\left(\underline{\tau},
\overline{\tau}\right),$
let
\begin{equation*}
  \underline{v}(v_0,z_0,\alpha)
  =\lim_{\tau\rightarrow\underline{\tau}(v_0,z_0,\alpha)}
  v(\tau;v_0,z_0,\alpha),
  \quad
  \overline{v}(v_0,z_0,\alpha)
  =\lim_{\tau\rightarrow\overline{\tau}(v_0,z_0,\alpha)}
  v(\tau;v_0,z_0,\alpha).
\end{equation*}
Then, we see that the trajectory
$L$ on
$\left(\underline{\tau},
\overline{\tau}\right)$
is also the solution curve
of the first-order differential equation
\begin{equation}\label{transformedSys:dz/dv}
  \frac{\mathrm{d}z}{\mathrm{d}v}
  =\frac{zf(v,\alpha)}{vf(v,\alpha)+z^2}
\end{equation}
on $\left(\underline{v},\overline{v}\right)$
with initial point $(v_0,z_0).$
Denote the corresponding solution by
\begin{equation*}
  z=\phi(v;v_0,z_0,\alpha),
\end{equation*}
or $z=\phi(v)$ for short,
on $\left(\underline{v},
\overline{v}\right)$
with initial condition $z_0=\phi(v_0;v_0,z_0,\alpha).$



\begin{rmk}
\normalfont{
In fact,
$\left(\underline{v},\overline{v}\right)$
is the maximal existence interval of
the solution of equation
\eqref{transformedSys:dz/dv}
with initial point
$(v_0,z_0).$
}
\end{rmk}

In the upper half plane,
$\mathrm{d}z/\mathrm{d}v$ changes signs
when the trajectory crosses the lines
$v=v_i^0$ for $i=1,2,3.$ Note that if
a trajectory $L$ passes through
a point $(v_0,z_0)\in\mathcal{D},$
$v(\tau)$ monotonically
depends on $\tau$ on
$(\underline{\tau},\overline{\tau}).$
There is some $\tau^*$ such that
the segment of
$L$ on
$\left(\tau^*,\overline{\tau}\right)$
lies entirely in some domain $\mathcal{D}_i,$
in which $\mathrm{d}z/\mathrm{d}v$ does not change sign.
Hence, $z(\tau)$ monotonically varies
when $\tau$ eventually approaches
$\overline{\tau}.$
It follows that
$\phi(v)$ monotonically varies
when $v$ eventually
approaches $\overline{v}.$
Similarly, $\phi(v)$ monotonically varies
when $v$ approaches $\underline{v}.$
Moreover, we have the following result.

\begin{lem}\label{lemma_phi<infty_if_v<infty}
Assume that $\alpha<-3/\sqrt[3]{4}$
and that the solution (or the trajectory)
$L(v_0,z_0,\alpha)$
of system \eqref{transformedSys(v,z)}
passes through a point
$(v_0,z_0)\in \mathcal{D}(\alpha).$
If $\overline{v}(v_0,z_0,\alpha)<+\infty$
(resp. $\underline{v}(v_0,z_0,\alpha)>-\infty$),
then
$$
\lim_{\tau\rightarrow\overline{\tau}(v_0,z_0,\alpha)}
z(\tau;v_0,z_0,\alpha)<+\infty
\left(\text{resp.}\;
\lim_{\tau\rightarrow\underline{\tau}(v_0,z_0,\alpha)}
z(\tau;v_0,z_0,\alpha)<+\infty
\right).
$$
\begin{proof}
We only consider the case for $\overline{v}$
and the analysis on the case for
$\underline{v}$
is analogous.

From the previous analysis, we see that
the limit of $z(\tau)$
at $\overline\tau$ exists or is equal to $+\infty.$
Assume that
$\lim_{\tau\rightarrow\overline{\tau}}
z(\tau)=+\infty.$
It follows that
$$
\lim_{v\rightarrow\overline{v}}
\phi(v)=+\infty,
$$
where $\phi(v)$ is
the solution of \eqref{transformedSys(v,z)} on
$\left(\underline{v},\overline{v}\right)$
with initial condition
$z_0=\phi(v_0).$
Since $\overline{v}<+\infty,$
there is a sufficiently large $M$ such that
for some $\epsilon_0>0$
\begin{equation*}
  \left|\frac{\mathrm{d}z}{\mathrm{d}v}\right|
=\left|\frac{zf(v,\alpha)}{vf(v,\alpha)+z^2}\right|
<1,
\quad
\forall (v,z)\in
\left(\overline{v}-\epsilon_0,
\overline{v}\right)
\times (M,+\infty),
\end{equation*}
and that for some $\epsilon_1\in(0,\epsilon_0)$
$$
\phi(v)>M,
\quad\forall v\in
\left(\overline{v}-\epsilon_1,
\overline{v}\right).
$$
Take some
$\tilde v\in\left(\overline{v}-\epsilon_1,
\overline{v}\right),$
then we have
\begin{equation*}
\begin{aligned}
  \lim_{v\rightarrow \overline{v}}\phi(v)
  =&
  \lim_{v\rightarrow \overline{v}}
  \int_{\tilde v}^v\mathrm{d}\phi
  +\phi(\tilde v)
  \\
  \le&
  \lim_{v\rightarrow \overline{v}}
  \int_{\tilde v}^v
  \left|\frac{\phi(v)f(v,\alpha)}
  {vf(v,\alpha)+\left(\phi(v)\right)^2}\right|
  dv+\phi(\tilde v)
  <
  \phi(\tilde v)+\epsilon_1
  <+\infty,
\end{aligned}
\end{equation*}
which is a contradiction.
\end{proof}
\end{lem}

Given a trajectory
$L$ passing a point
$(v_0,z_0)\in\mathcal{D},$
we have two claims:
\begin{enumerate}[(i)]
  \item
  $L$ approaches the curve $\mathcal{K}^+$
  or the equilibrium $P_0$ or $P_3$
  as $\tau\rightarrow\overline{\tau}$
  if $\overline{v}<+\infty,$
  and $\lim_{\tau\rightarrow+\infty}
  z(\tau)=+\infty$
  if $\overline{v}=+\infty;$
  \item
  $L$ approaches the curve $\mathcal{K}^-$
  or the equilibrium $P_1$ or $P_2$
  as $\tau\rightarrow\underline{\tau}$
  if $\underline{v}>-\infty,$
  and $\lim_{\tau\rightarrow-\infty}
  z(\tau)=+\infty$
  if $\underline{v}=-\infty.$

\end{enumerate}
For the case $\overline{v}<+\infty,$
since from
Lemma \ref{lemma_phi<infty_if_v<infty}
$\lim_{\tau\rightarrow\overline{\tau}}
z(\tau)$ always exists,
$L$ must approach some point as
$\tau\rightarrow\overline{\tau}.$
Precisely, if
$\overline{v}<+\infty$
and $\overline{\tau}<+\infty,$
then the regular point
$
\left(
v(\overline{\tau}),
z(\overline{\tau})
\right)
$
is obviously on $\mathcal{K}^+;$
if
$\overline{v}<+\infty$
and $\overline{\tau}=+\infty,$
then the limit point is
the equilibrium, i.e., $P_0$ or $P_3.$
For the case
$\overline{v}=+\infty$
which implies that
$\overline{\tau}=+\infty,$
there is a $\tau^*$ such that
the segment of $L$ on
$(\tau^*,+\infty)$
lies in
$\mathcal{D}_4$
in which
$\mathrm{d}z/\mathrm{d}\tau
\ge f(\tilde v,\alpha)z$
for any $v\ge \tilde v\gg0,$
and hence,
$$
z(\tau)\ge
z(\tau^*)\exp\left[f(\tilde v,\alpha)(\tau-\tau^*)\right]
\rightarrow+\infty
\quad
\text{as}
\quad
\tau\rightarrow+\infty.
$$
The discussion on $\underline{v}$ is analogous.

\subsection{Behavior of the separatrices of
$P_1(\alpha)$ and $P_3(\alpha)$}

Proposition \ref{proposition_equilibria_of_sys(v,z)}
implies that the unstable separatrix of $P_1(\alpha)$
in the upper half plane, denoted by $L^-(\alpha),$
passes through $\mathcal{D}_2(\alpha),$
and the stable separatrix of $P_3(\alpha)$ in the upper half plane,
denoted by $L^+(\alpha),$ passes through
$\mathcal{D}_3(\alpha).$
Choose two points
$(v_1,z_1)\in\mathcal{D}_2(\alpha)\cap L^-(\alpha)$
and
$(v_2,z_2)\in\mathcal{D}_3(\alpha)\cap L^+(\alpha),$
then we have
\begin{equation*}
  L^-(\alpha)=L(v_1,z_1,\alpha),\quad
  L^+(\alpha)=L(v_2,z_2,\alpha).
\end{equation*}
Let
\begin{equation*}
  v^-(\alpha)=\overline{v}(v_1,z_1,\alpha),\quad
  v^+(\alpha)=\underline{v}(v_2,z_2,\alpha),
\end{equation*}
and
\begin{equation*}
  \phi^-(v;\alpha)=\phi(v;v_1,z_1,\alpha),
  \quad
  \phi^+(v;\alpha)=\phi(v;v_2,z_2,\alpha).
\end{equation*}
Note that these notations are independent of the
choices of $(v_1,z_1)$ and $(v_2,z_2).$
Since $L^-(\alpha)$ negatively approaches $P_1(\alpha),$
we have $\underline{v}(v_1,z_1,\alpha)=v_1^0(\alpha),$
and hence, $\phi^-(v;\alpha)$ is defined on
$v\in\left(v_1^0(\alpha),v^-(\alpha)\right).$
Similarly, $\phi^+(v;\alpha)$ is defined on
$v\in\left(v^+(\alpha),v_3^0(\alpha)\right).$

In terms of
Proposition \ref{proposition_equilibria_of_sys(v,z)}
and the claims after Lemma \ref{lemma_phi<infty_if_v<infty},
we can obtain the
the possible behavior of $L^-(\alpha)$ and $L^+(\alpha),$
which are entirely related to the possible values of
$v^-(\alpha)$ and $v^+(\alpha).$

\begin{prop}\label{proposition_range_of_v+-}
Assume that $\alpha<-3/\sqrt[3]{4}.$
$v^-(\alpha)$ and $v^+(\alpha)$ admit one of the following cases:
\begin{enumerate}[(i)]
  \item
  $v^-(\alpha)=0$ and
  $v^+(\alpha)=-\infty,$
  see Figure \ref{fig:sys(v,z)_alpha=2};

  \item
  $0<v^-(\alpha)<v_3^0(\alpha)$ and
  $v^+(\alpha)=-\infty,$
  see Figure \ref{fig:sys(v,z)_alpha=2.15};

  \item
  $v^-(\alpha)=v_3^0(\alpha)$
  and
  $v^+(\alpha)=v_1^0(\alpha),$
  see Figure \ref{fig:sys(v,z)_alpha=2.198};

  \item
  $v^-(\alpha)=+\infty$
  and
  $v_1^0<v^+(\alpha)<v_2^0(\alpha),$
  see Figure \ref{fig:sys(v,z)_alpha=2.3}.
\end{enumerate}

\end{prop}

\begin{figure}[htbp]
\centering
\subfigure[]{
\label{fig:sys(v,z)_alpha=2}
\includegraphics[width=0.45\linewidth]{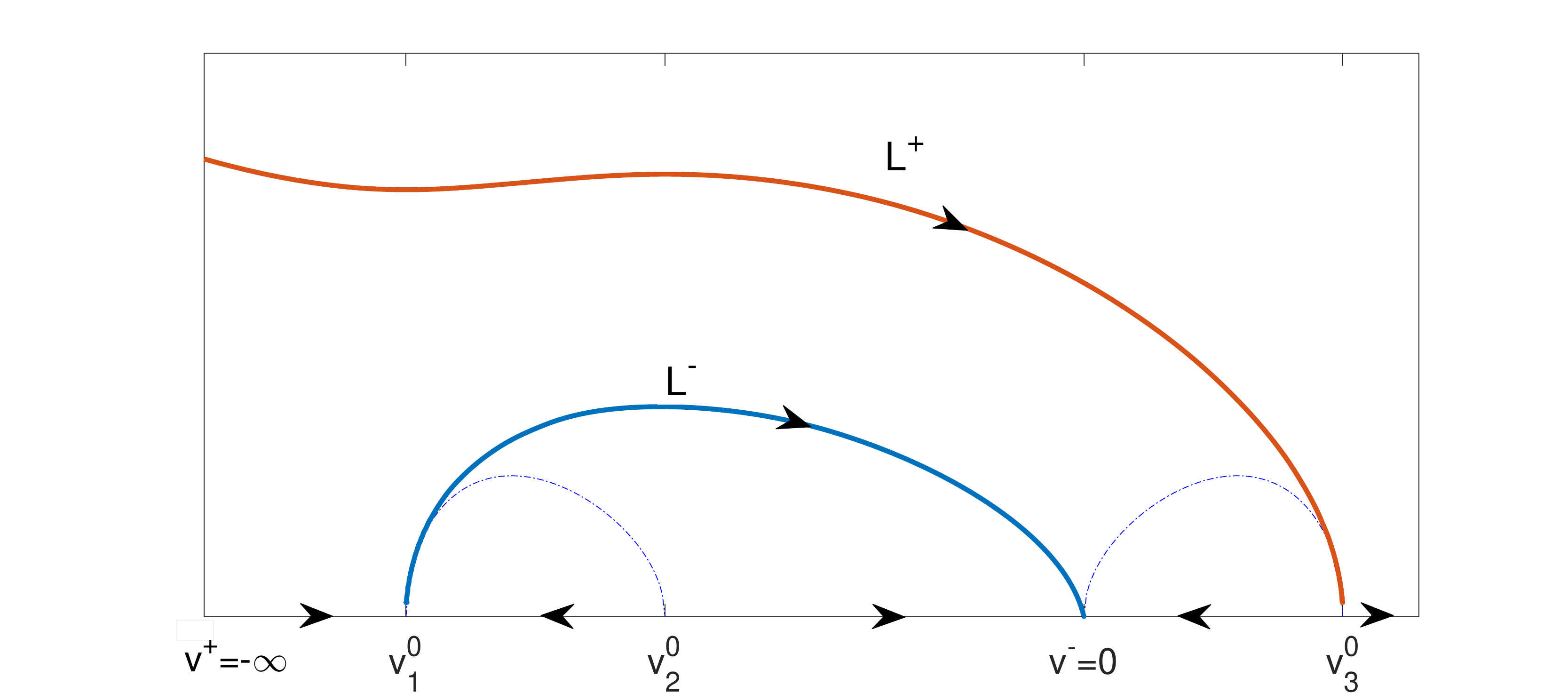}
}
\subfigure[]{
\label{fig:sys(v,z)_alpha=2.15}
\includegraphics[width=0.45\linewidth]{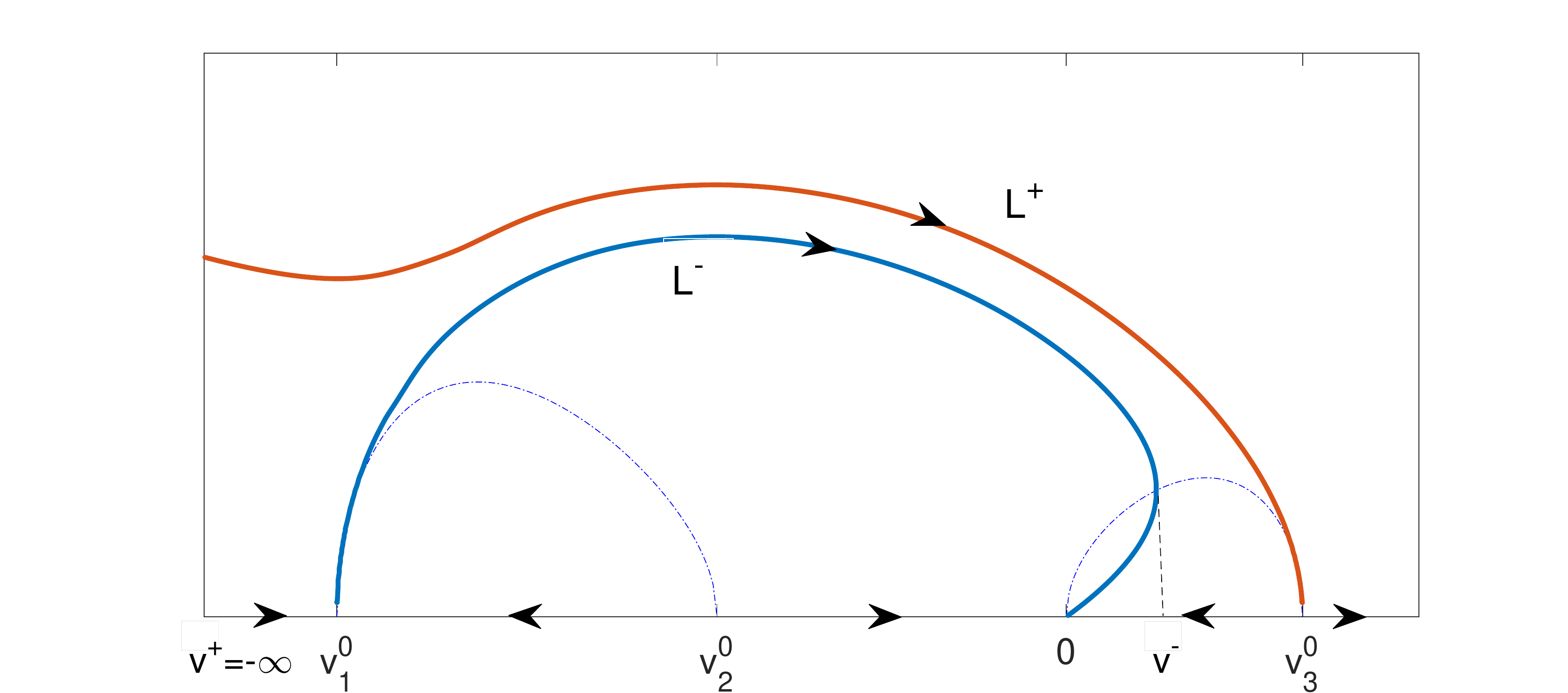}
}
\\
\subfigure[]{
\label{fig:sys(v,z)_alpha=2.198}
\includegraphics[width=0.45\linewidth]{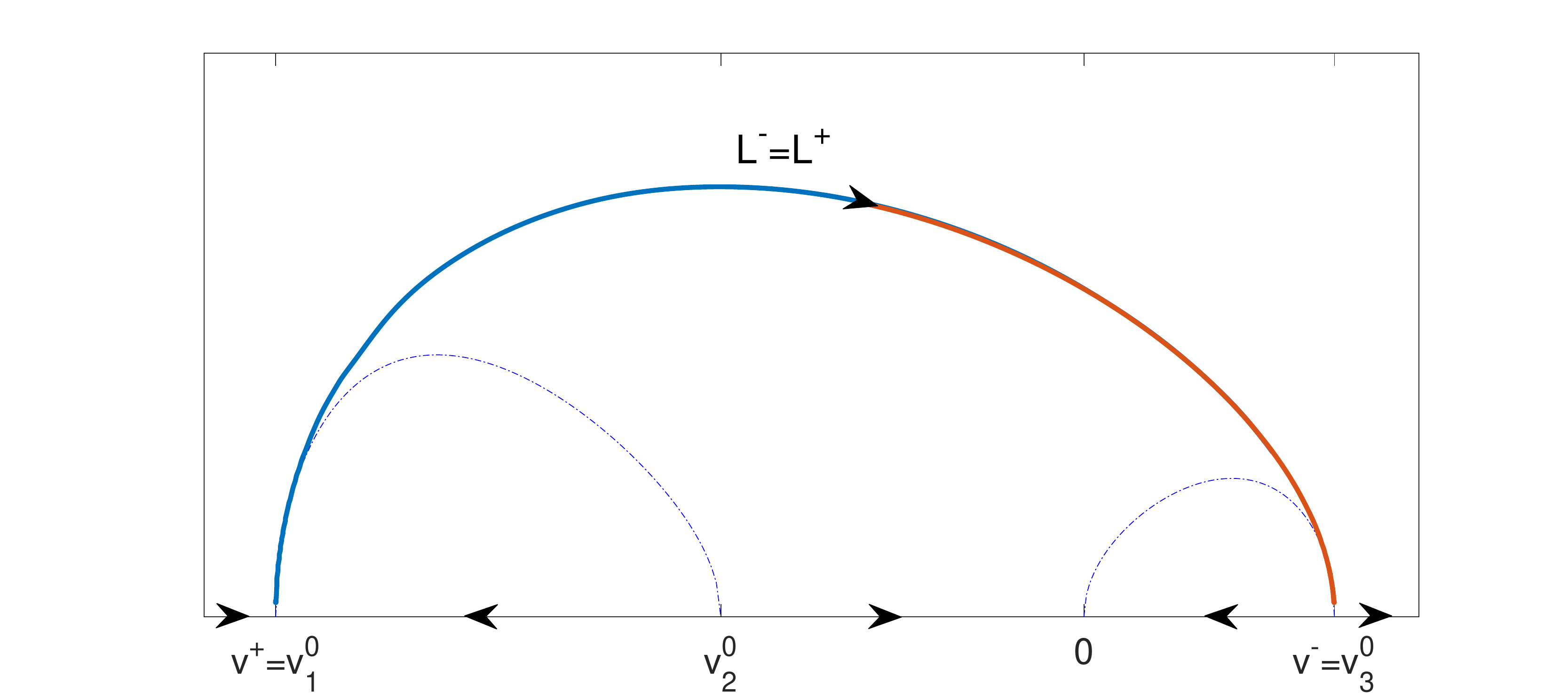}
}
\subfigure[]{
\label{fig:sys(v,z)_alpha=2.3}
\includegraphics[width=0.45\linewidth]{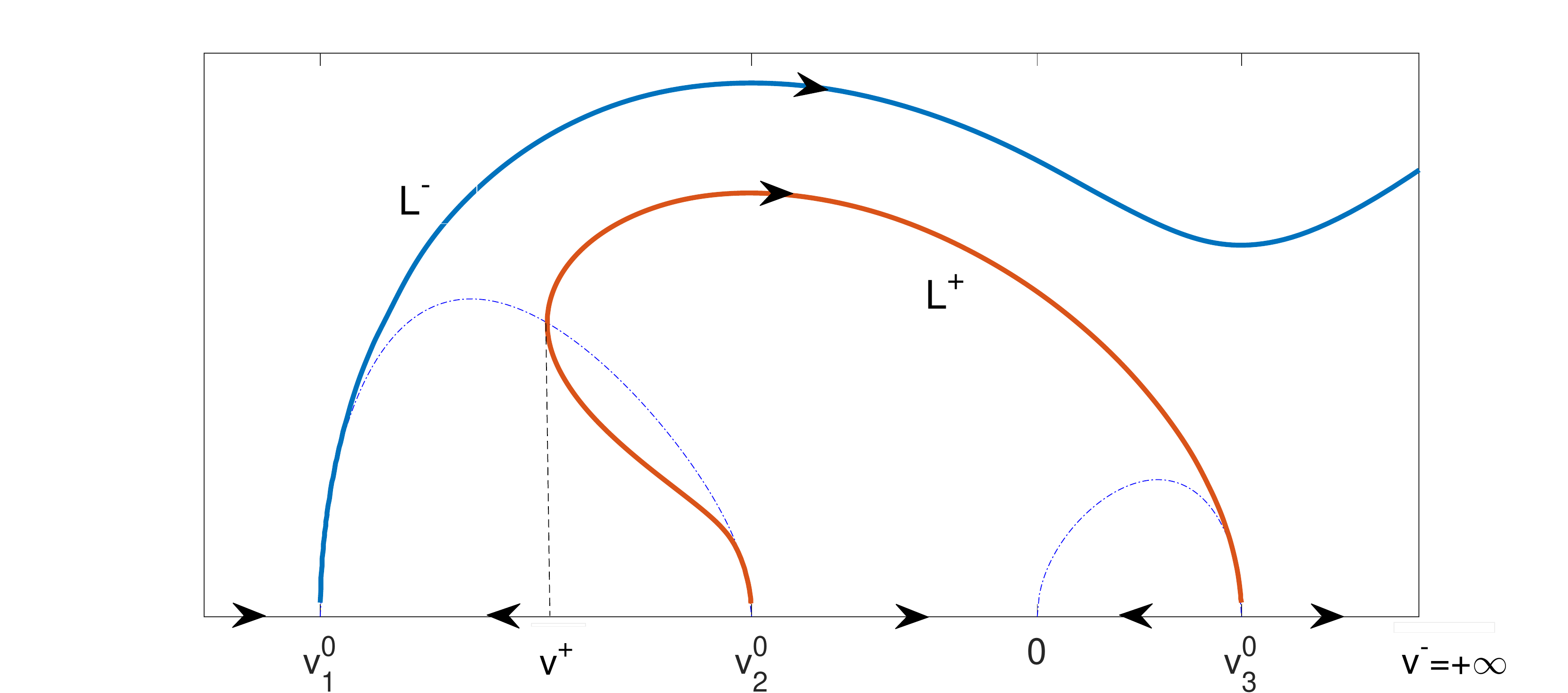}
}
\caption{Possible cases for the two separatrices $L^-$
and $L^+$ of the saddles $P_1$ and $P_3$
with $\alpha<-3/\sqrt[3]4,$ respectively.}
\label{fig:sys(v,z)}
\end{figure}


\begin{rmk}\label{remark_case_for_phi+-}
\normalfont{
From Proposition \ref{proposition_range_of_v+-},
one can see that
$\lim_{v\rightarrow0}\phi^-(v;\alpha)$ and
$\lim_{v\rightarrow0}\phi^+(v;\alpha)$
always exist, the former of
which is written as $\phi^-(0;\alpha)$ even though
$v^-(\alpha)=0,$
i.e., the case shown in Figure \ref{fig:sys(v,z)_alpha=2},
while the latter of which is always larger than 0.
Then the four cases listed in Proposition
\ref{proposition_range_of_v+-} correspond to
$\phi^-(0;\alpha)=0,$
$0<\phi^-(0;\alpha)<\phi^+(0;\alpha),$
$\phi^-(0;\alpha)=\phi^+(0;\alpha)$
and
$\phi^-(0;\alpha)>\phi^+(0;\alpha),$
respectively.

}
\end{rmk}

It is well known that any trajectory of
system \eqref{transformedSys(v,z)}
with fixed initial point continuously
depends on $\alpha$ in some suitable range.
That is to say,
$v(\tau;v_0,z_0,\alpha)$ and
$z(\tau;v_0,z_0,\alpha)$ continuously
depend on $\alpha$ in some suitable range
at any finite $\tau.$
However, this property is not much easy to
describe and confirm for the separatrices
$L^-(\alpha)$ and $L^+(\alpha),$ because
we only know that their limit points
$P_1(\alpha)$
and
$P_3(\alpha)$
continuously depend on $\alpha.$
Therefore, we first give a proof of the
continuity with respect to $\alpha$
of these two separatrices.

\begin{lem}\label{lemma_alpha_continuity_of_phi-+}
Assume that $\alpha_0<-3/\sqrt[3]{4}.$
For any $v_0\in\left(v_1^0(\alpha_0),0\right],$
$\phi^-(v_0;\alpha)$ continuously depends on
$\alpha$ at $\alpha_0.$
For any
$v_0\in\left[0,v_3^0(\alpha_0)\right),$
$\phi^+(v_0;\alpha)$ continuously depends on
$\alpha$ at $\alpha_0.$

\begin{proof}
We only consider the case for $\phi^-,$
and the analysis on the case for $\phi^+$ is analogous.

First consider the case $v^-(\alpha_0)>0,$
for which $\phi^-(v;\alpha_0)>0$ on
$ \left(v_1^0(\alpha_0),0\right]$
(see Figure \ref{fig:sys(v,z)_alpha=2.15}
to \ref{fig:sys(v,z)_alpha=2.3}).
For any $v_0\in\left(v_1^0(\alpha_0),0\right],$
since $v_1^0(\alpha)$ is continuous
with respect to $\alpha<-3/\sqrt[3]{4},$
there must be a $\delta_0>0$
such that $v_1^0(\alpha)<v_0$
on $U_{\delta_0}(\alpha_0),$
where $U_{\delta}(P)$
denotes the $\delta$-neighborhood
of a point $P.$
Hence, $\phi^-(v_0;{\alpha})$
makes sense for any
$\alpha\in U_{\delta_0}(\alpha_0).$

Let $z_0=\phi^-(v_0;\alpha_0).$
Since the point $(v_0,z_0)$ belongs to
the open set $\mathcal{D}(\alpha_0),$
in which $\mathrm{d}v/\mathrm{d}\tau>0,$
for any sufficiently small $\varepsilon>0$
we have
$U_{2\varepsilon}(v_0,z_0)\subset\mathcal{D}(\alpha_0).$
Consider the trajectory starts from
$(v_0,z_0+\varepsilon),$
i.e., $L(v_0,z_0+\varepsilon,\alpha_0).$
From Lemma \ref{lemma_phi<infty_if_v<infty}
and
Proposition \ref{proposition_equilibria_of_sys(v,z)}(ii),
we see that
$L(v_0,z_0+\varepsilon,\alpha_0)$
negatively goes into
the inner of $\mathcal{D}_1(\alpha_0).$
That is to say, there is a $\tau_1<0$ such that
$$
v(\tau_1;v_0,z_0+\varepsilon,\alpha_0)
<v_1^0(\alpha_0),\quad
z(\tau_1;v_0,z_0+\varepsilon,\alpha_0)>0.
$$
Since $v(\tau_1;v_0,z_0+\varepsilon,\alpha),
z(\tau_1;v_0,z_0+\varepsilon,\alpha)$
and $v_1^0(\alpha)$
are continuous respect to $\alpha,$
there is a $\delta_1>0$ such that
\begin{equation}\label{lemma_alpha_continuity_of_phi-+_ineq1}
  \left(
v(\tau_1;v_0,z_0+\varepsilon,\alpha),
z(\tau_1;v_0,z_0+\varepsilon,\alpha)
\right)
\in\mathcal{D}_1(\alpha),\quad
\forall \alpha\in U_{\delta_1}(\alpha_0).
\end{equation}
On the other hand,
consider the trajectory starts from
$(v_0,z_0-\varepsilon),$
i.e., $L(v_0,z_0-\varepsilon,\alpha_0).$
From Proposition \ref{proposition_equilibria_of_sys(v,z)}(ii)
and (iii),
we see that $L(v_0,z_0-\varepsilon,\alpha_0)$
negatively goes into $\mathcal{A}_1(\alpha_0)$
and eventually approaches $P_2(\alpha_0).$
That is to say, there is a $\tau_2<0$ such that
$$
\begin{aligned}
&v_1^0(\alpha_0)<
v(\tau_2;v_0,z_0-\varepsilon,\alpha_0)
<v_2^0(\alpha_0),
\\
&0<z(\tau_2;v_0,z_0-\varepsilon,\alpha_0)
<\sqrt{-v(\tau_2;v_0,z_0-\varepsilon,\alpha_0)
f\left(v(\tau_2;v_0,z_0-\varepsilon,\alpha_0),\alpha_0\right)}.
\end{aligned}
$$
Since $v(\tau_2;v_0,z_0-\varepsilon,\alpha),
z(\tau_2;v_0,z_0-\varepsilon,\alpha),$
$v_1^0(\alpha),v_2^0(\alpha)$
and $\mathcal{K}^-(\alpha)$ are continuous
respect to $\alpha,$ there is a $\delta_2>0$
such that
\begin{equation}\label{lemma_alpha_continuity_of_phi-+_ineq2}
  \left(
v(\tau_2;v_0,z_0-\varepsilon,\alpha),
z(\tau_2;v_0,z_0-\varepsilon,\alpha)
\right)
\in\mathcal{A}_1(\alpha),\quad
\forall \alpha\in U_{\delta_2}(\alpha_0).
\end{equation}

Now take
$\delta_3=\min\{\delta_0,\delta_1,\delta_2\}.$
In terms of
\eqref{lemma_alpha_continuity_of_phi-+_ineq1},
\eqref{lemma_alpha_continuity_of_phi-+_ineq2}
and Proposition
\ref{proposition_equilibria_of_sys(v,z)}(ii),
we see that $L^-(\alpha)$ is always bounded by
$L(v_0,z_0+\varepsilon,\alpha)$ and
$L(v_0,z_0-\varepsilon,\alpha)$
for any $\alpha\in U_{\delta_3}(\alpha_0).$
Hence, we have
$$
z_0-\varepsilon<\phi^-(v_0;\alpha)
<z_0+\varepsilon,\quad
\forall \alpha\in U_{\delta_3}(\alpha_0),
$$
which completes the proof for the case
$v^-(\alpha_0)>0.$

If $v^-(\alpha_0)=0,$
then $\phi^-(0;\alpha_0)=0$
(see Figure \ref{fig:sys(v,z)_alpha=2}).
We only need to consider $v_0=0.$
For any $\varepsilon>0,$
there is a $\tilde v(<0)$ sufficiently closed
to 0 such that
\begin{equation*}
  \phi^-(\tilde v;\alpha_0)<\frac{\varepsilon}{2}.
\end{equation*}
From the first case,
there is $\delta>0$ such that
\begin{equation*}
  \left|
  \phi^-(\tilde v;\alpha)-\phi^-(\tilde v;\alpha_0)
  \right|
  <\frac{\varepsilon}{2},\quad
  \forall\alpha\in U_{\delta}(\alpha_0).
\end{equation*}
Recall that
$\mathrm{d}z/\mathrm{d}\tau<0$
and
$\mathrm{d}v/\mathrm{d}\tau>0$
in the inner of
$\mathcal{D}_3(\alpha).$
It follows that
$\mathrm{d}\phi^-(v;\alpha)/\mathrm{d}v<0$
on $\left(v_2^0(\alpha),0\right).$
Then, we have
\begin{equation*}
  \phi^-(0;\alpha)
  <\phi^-(\tilde v;\alpha)
  \le\left|
  \phi^-(\tilde v;\alpha)-\phi^-(\tilde v;\alpha_0)
  \right|
  +\phi^-(\tilde v;\alpha)
  <\varepsilon,\quad
  \forall \alpha\in U_{\delta}(\alpha_0).
\end{equation*}
The proof is completed.
\end{proof}
\end{lem}

%

\begin{rmk}
\normalfont{
Given a smooth curve $\mathcal{C},$
it is known that the intersection of
a trajectory $L(v,z,\alpha)$
and $\mathcal{C}$ varies
continuously on $\mathcal{C}$
when $v,z$ and $\alpha$ vary
in some suitable ranges.
The continuity of $L^-(\alpha)$
 with respect to
$\alpha$ is described in the same manner
and it can be proved
by applying
Lemma \ref{lemma_alpha_continuity_of_phi-+}
on
$
L^-(\alpha)=L
\left(v_0,\phi^-(v_0;\alpha),\alpha\right)
$
with some fixed
$v_0\in\left(v_1^0(\alpha_0),0\right).$
Similar statements hold for $L^+(\alpha).$
However, we only need to consider the continuity
of $\phi^-(0;\alpha)$ and $\phi^+(0;\alpha)$
with respect to $\alpha$ in this paper.
}
\end{rmk}

For the main purposes of this section,
we are going to show
(in Lemma \ref{lemma_L_monotonically_varies})
that
$\phi^-(0;\alpha)$ (if larger than 0) and $\phi^+(0;\alpha)$
vary monotonically as $\alpha$ varies,
and (in Lemma \ref{lemma_L-_above_L+})
that they behave as Figure \ref{fig:sys(v,z)_alpha=2.3}
when $\alpha=-3/\sqrt[3]{2}.$
Before this, let's introduce a useful
theorem which is a modified version for the
\emph{basic comparison theorem} in \cite{mcnabb1986comparison}.
\begin{thm}[Comparison theorem]\label{thm_comparison}
Let $\xi(x)$ and $\eta(x)$ be continuously differentiable
functions on $(x_0,x_1),$
and $g(x,y)$ be a continuous function on $D\subset\mathbb{R}^2,$
where $D$ contains the set
$
\left\{(x,\xi(x)) | x_0< x<x_1 \right\}
\cup
\left\{(x,\eta(x)) | x_0< x<x_1 \right\}.
$
Assume that
  \begin{equation}\label{thm_comparison_cond1}
  \frac{\mathrm{d}\xi}{\mathrm{d}x}-g\left(x,\xi(x)\right)
  <
  \frac{\mathrm{d}\eta}{\mathrm{d}x}-g\left(x,\eta(x)\right),
  \quad
  \forall x\in(x_0,x_1),
  \end{equation}
  and either
  \begin{equation}\label{thm_comparison_cond2}
  \limsup_{x\rightarrow x_0}
  \left(\xi(x)-\eta(x)\right)<0
  \end{equation}
  or
  \begin{equation}\label{thm_comparison_cond3}
  \limsup_{x\rightarrow x_0}
  \left(\xi(x)-\eta(x)\right)=0,
  \quad
  \limsup_{x\rightarrow x_0}\left(
  \frac{\mathrm{d}\xi}{\mathrm{d}x}
  -\frac{\mathrm{d}\eta}{\mathrm{d}x}
  \right)<0.
  \end{equation}
Then $\xi(x)<\eta(x)$ on $(x_0,x_1).$
Similar statements hold for $x_1<x_0.$
\begin{proof}
Assume that $\xi(x)\ge\eta(x)$ for some $x\in(x_0,x_1).$
We can take
$$
p^*=\inf\{x\in(x_0,x_1) \;|\; \xi(x)\ge\eta(x)\}.
$$
If $p^*>x_0,$ then from the continuity of
$\xi$ and $\eta$ on $(x_0,x_1),$ we have
$\xi(p^*)=\eta(p^*).$
However, from assumption \eqref{thm_comparison_cond1},
we have
$$
\left.\frac{\mathrm{d}\xi}{\mathrm{d}x}\right|_{x=p^*}
-g\left(p^*,\xi(p^*)\right)
  <
\left.\frac{\mathrm{d}\eta}{\mathrm{d}x}\right|_{x=p^*}
-g\left(p^*,\eta(p^*)\right),
$$
and hence,
$$
\left.\frac{\mathrm{d}\xi}{\mathrm{d}x}\right|_{x=p^*}
<
\left.\frac{\mathrm{d}\eta}{\mathrm{d}x}\right|_{x=p^*}.
$$
Thus, there is a $x\in(x_0,p^*)$ such that
$\xi(x)>\eta(x),$
which contradicts to the definition of $p^*.$

Then, we have $p^*=x_0.$ From the definition of $p^*,$
there is a sequence $\{p_n\}$ on $(x_0,x_1)$
such that
$$
p_{n+1}<p_{n},
\quad
\lim_{n\rightarrow\infty}p_n=x_0,
\quad
\zeta(p_n)\ge0,
$$
where $\zeta(x)=\xi(x)-\eta(x).$
It follows that
$$
\limsup_{x\rightarrow x_0}\zeta(x)
\ge
\limsup_{n\rightarrow \infty}\zeta(p_n)
\ge
\liminf_{n\rightarrow \infty}\zeta(p_n)
\ge0,
$$
which contradicts to condition
\eqref{thm_comparison_cond2}.
Thus, we must have
$\limsup_{x\rightarrow x_0}\zeta(x)=0,$
and hence,
\begin{equation}\label{thm_comparison_limitzetapn=0}
\lim_{n\rightarrow \infty}\zeta(p_n)=0.
\end{equation}
On the other hand, the latter part of condition \eqref{thm_comparison_cond3}
implies that there is a
sufficiently small $\epsilon>0$ such that
$$
\frac{\mathrm{d}\zeta}{\mathrm{d}x}<0,
\quad
\forall x\in(x_0,x_0+\epsilon).
$$
Then, we have
$$
\zeta(p_{n+1})-\zeta(p_n)=(p_{n+1}-p_n)
\left.\frac{\mathrm{d}\zeta}{\mathrm{d}x}\right|_{x=q_n}
>0
$$
for some $q_n\in[p_{n+1},p_n]$ and any sufficiently large $n$
such that $p_n<x_0+\epsilon.$
Therefore, when $p_n<x_0+\epsilon,$ we have
$$
0\le\zeta(p_n)<\zeta(p_{n+1})<\zeta(p_{n+2})<\cdots,
$$
which contradicts to \eqref{thm_comparison_limitzetapn=0}.
The disproof of the theorem is completed.

Similar statements hold for $x_1<x_0$ by considering
$-x\in(-x_0,-x_1).$
\end{proof}
\end{thm}

\begin{rmk}\label{remark_comparison}
\normalfont{
In particular, if $\xi(x)$ and $\eta(x)$
are solutions of two first-order differential equations
$\mathrm{d}y/\mathrm{d}x=g_1(x,y)$
and
$\mathrm{d}y/\mathrm{d}x=g_2(x,y),$
respectively, and if
$g_1(x,\eta(x))$ (resp. $g_2(x,\xi(x))$)
makes sense on $(x_0,x_1),$
then condition \eqref{thm_comparison_cond1}
in Theorem \ref{thm_comparison} reduces to
$$
g_1(x,\eta(x))<g_2(x,\eta(x)) \quad
\left(\text{resp.}\;\;
g_1(x,\xi(x))<g_2(x,\xi(x))
\right)
$$
by taking $g=g_1$ (resp. $g=g_2$).
}
\end{rmk}


\begin{lem}\label{lemma_L_monotonically_varies}
Assume that $b<a<-3/\sqrt[3]{4}.$
For any
$v\in\left(v_1^0(a),0\right],$
we have
\begin{equation*}
    \phi^-(v;a)
    \le\phi^-(v;b),
\end{equation*}
where the equal sign holds if and only if
$v=v^-(b)=0$
(which is equivalent to
 $\phi^-(0;b)=0$
from Remark
\ref{remark_case_for_phi+-} ).
For any $v\in\left[0,v_3^0(b)\right),$
we have
\begin{equation*}
    \phi^+(v;a)>\phi^+(v;b).
\end{equation*}

\begin{proof}
We only consider $\phi^{-}$
and the analysis on $\phi^{+}$ is analogous.

First, consider the case $v\in\left(v_1^0(a),0\right).$
Recall that $v_1^0(a)>v_1^0(b)$ from Lemma \ref{lemma_features_of_vmui},
which implies that
\begin{equation}\label{lemma_L_monotonically_varies_eq1}
\lim_{v\rightarrow v_1^0(a)}\phi^-(v;a)=0
<\phi^-\left(v_1^0(a);b\right).
\end{equation}
Note that when $z>0$ and $vf(v,\alpha)+z^2\neq0,$ we have
$$
\frac{\partial }{\partial \alpha} \frac{zf(v,\alpha)}{vf(v,\alpha)+z^2}
=
-\frac{v^2z^3}{\left(vf(v,\alpha)+z^2\right)^2}<0,
$$
and that
$
0<\left(\phi^-(v;b)\right)^2+vf(v,b)
\le\left(\phi^-(v;b)\right)^2+vf(v,\alpha)
$
for any
$\alpha\in[b,a],$
which implies that
$\left(v,\phi^-(v;b)\right)
\in\mathcal{D}(\alpha)$
for any $\alpha\in[b,a].$
It follows that
\begin{equation}\label{lemma_L_monotonically_varies_eq2}
  \left.
  \frac{zf(v,a)}{vf(v,a)+z^2}
  \right|_{z=\phi^-(v;b)}
  <
  \left.
  \frac{zf(v,b)}{vf(v,b)+z^2}
  \right|_{z=\phi^-(v;b)}.
\end{equation}
In terms of
\eqref{lemma_L_monotonically_varies_eq1},
\eqref{lemma_L_monotonically_varies_eq2}
and applying Theorem \ref{thm_comparison} and
Remark \ref{remark_comparison} on $(v_1^0(a),0),$
we have
$\phi^-(v;a)<\phi^-(v;b)$
for any $v\in\left(v_1^0(a),0\right).$

Next, consider the case $v=0$ but $v^-(b)>0$
which implies that $\phi^-(0;b)>0$
from Remark
\ref{remark_case_for_phi+-}.
From the first case in this proof, we have
$$
\phi^-(0;a)
=\lim_{v\rightarrow0,v<0}\phi^-(0;a)
\le\lim_{v\rightarrow0,v<0}\phi^-(0;b)
=\phi^-(0;b).
$$
Assume that $\phi^-(0;a)=\phi^-(0;b).$
Then
$\left(0,\phi^-(0;a)\right)$
is a regular point on both $L^-(a)$ and $L^-(b).$
Let $\tilde D$ be the open domain bounded by
$L^{-}(a),L^{-}(b)$ and the line segment
$\overline{P_1(b)P_1(a)}.$
Let $\tilde L$ be a trajectory of
system \eqref{transformedSys(v,z)} with $\alpha=a$
passing through $\tilde D.$
Since both $L^{-}(a)$ and $\overline{P_1(b)P_1(a)}$
are trajectories of system \eqref{transformedSys(v,z)} with $\alpha=a,$
$\tilde L$ must cross $L^{-}(b)$ from $\tilde D$
to $\mathbb{R}^2-\tilde D$ at a point
$\tilde P\in L^{-}(b)\cap\{v<0\}.$
Together with
$\left.\mathrm{d}v/\mathrm{d}\tau\right|_{\tilde P} >0,$
i.e.,
$\tilde P\in\mathcal{D}(\alpha),$
for any $\alpha\in[b,a],$ we have
$$
\left.\frac{\mathrm{d}z}{\mathrm{d}v}\right|_{\tilde P,\alpha=a}
\ge\left.\frac{\mathrm{d}z}{\mathrm{d}v}\right|_{\tilde P,\alpha=b}
\quad\text{but}\quad
\left.
\frac{\partial }{\partial \alpha} \frac{\mathrm{d}z}{\mathrm{d}v}
\right|_{\tilde P}<0,
\;\forall\alpha\in[b,a],
$$
which is a contradiction.

At last, from the aforementioned analysis, the equal
sign only occurs in the case
$\phi^-(0;a)=\phi^-(0;b)=0,$
which is equivalent to
$v^-(a)=v^-(b)=0$
from Remark
\ref{remark_case_for_phi+-}.

\end{proof}
\end{lem}

Next we show that when $\alpha=-3/\sqrt[3]{2}$
the two separatrices behave as the case
shown in Figure \ref{fig:sys(v,z)_alpha=2.3}.
\begin{lem}\label{lemma_L-_above_L+}
$
\phi^-\left(0;-3/\sqrt[3]{2}\right)
>\phi^+\left(0;-3/\sqrt[3]{2}\right).
$

\begin{proof}
We estimate
$\phi^-(0;\alpha)$
and
$\phi^+(0;\alpha)$
by constructing two curves such that
one lies below $L^{-}(\alpha)$
and the other above $L^{+}(\alpha).$
For convenience, if $\alpha$ is arbitrary,
$\phi^-(v;\alpha),$
$\phi^+(v;\alpha)$
and $f(v,\alpha)$
are abbreviated to
$\phi^-(v),\phi^+(v)$
and $f(v),$ respectively.

Let $\alpha<-3/\sqrt[3]{4}$ and
$0<\mu<f\left(2\alpha/3\right)=-4\alpha^3/27-1.$
Consider the following piecewise differential equation
\begin{equation*}\label{piecewise diffeq-}
\frac{\mathrm{d}z}{\mathrm{d}v}
=G^{\mu}(v,z)\triangleq
\left\{
       \begin{aligned}
       &\frac{\mu z}{\mu v+z^2},  &v\in\left[v^{\mu}_1,v^{\mu}_2\right),\\
       &0,                        &v\in\left[v^{\mu}_2,v^{0}_2\right),\\
       &\frac{z}{v-z^2},                      &v\in\left[v^{0}_2,0\right],
       \end{aligned}
\right.
\end{equation*}
where $v^{\mu}_i$ is defined by \eqref{def_of_vmui}.
Take the point
$\left(v^{\mu}_1,z_1^{\mu}\right)$
between the curve $\mu v+z^2=0$ and $L^{-},$
i.e.,
\begin{equation}
\label{lemma_L-_above_L+_initial_cond1_for_phi-}
  \sqrt{-\mu v^{\mu}_1}<z_1^{\mu}<\phi^-\left(v^{\mu}_1\right).
\end{equation}
Let $z=\varphi^{\mu}(v)$ be the solution of equation
$\mathrm{d}z/\mathrm{d}v=G^{\mu}(v,z)$
with initial point
$\left(v^{\mu}_1,z_1^{\mu}\right),$
i.e.,
$z_1^\mu=\varphi^\mu\left(v^{\mu}_1\right).$
It is not hard to show that
$\varphi^{\mu}(v)$ exists on $\left[v^{\mu}_1,0\right].$
In fact, it can be explicitly formulated by
\begin{equation*}\label{explicit_sol_of_varphimu}
  \varphi^{\mu}(v)=
  \left\{
  \begin{aligned}
        &\frac{1}{2}\left[
        z_1^{\mu}-\mu\frac{v_1^{\mu}}{z_1^{\mu}}
        +\sqrt{\left(z_1^{\mu}+\mu\frac{v_1^{\mu}}{z_1^{\mu}}\right)^2
               +4\mu\left(v-v_1^{\mu}\right)}
        \right],
        &v\in\left[v^{\mu}_1,v^{\mu}_2\right),\\
        &z_2^{\mu},
        &v\in\left[v^{\mu}_2,v^{0}_2\right),\\
        &\frac{1}{2}\left[
        z_2^{\mu}+\frac{v_2^{0}}{z_2^{\mu}}
        +\sqrt{\left(z_2^{\mu}+\frac{v_2^{0}}{z_2^{\mu}}\right)^2-4v}
        \right],
        &v\in\left[v^{0}_2,0\right],
  \end{aligned}
  \right.
\end{equation*}
where
\begin{equation*}
\begin{aligned}
  z_2^{\mu}=&
  \frac{1}{2}\left[
        z_1^{\mu}-\mu\frac{v_1^{\mu}}{z_1^{\mu}}
        +\sqrt{\left(z_1^{\mu}+\mu\frac{v_1^{\mu}}{z_1^{\mu}}\right)^2
               +4\mu\left(v_2^{\mu}-v_1^{\mu}\right)}
        \right] \\
  \ge&
  \sqrt{\mu}\left(\sqrt{-v_1^{\mu}}
  +\sqrt{v_2^{\mu}-v_1^{\mu}}\right).
\end{aligned}
\end{equation*}
Notice
\begin{equation*}
f(v)>\left\{
  \begin{aligned}
  &\mu, &v\in\left(v^{\mu}_1,v^{\mu}_2\right),\\
  &0, &v\in\left[v^{\mu}_2,v^{0}_2\right),\\
  &-1, &v\in\left[v^{0}_2,0\right),
  \end{aligned}
\right.
\end{equation*}
$f(v)<0$ on $\left(v^{0}_2,0\right)$ and
\begin{equation*}
  \begin{aligned}
  &0<v+\frac{\left(\phi^-(v)\right)^2}{f(v)}
  <v+\frac{\left(\phi^-(v)\right)^2}{\mu},\quad
  &\forall v\in\left(v_1^\mu,v_2^\mu\right),
  \\
  &v+\frac{\left(\phi^-(v)\right)^2}{f(v)}
  <v-\left(\phi^-(v)\right)^2<0,\quad
  &\forall v\in\left(v_2^0,0\right).
  \end{aligned}
\end{equation*}
It follows that
\begin{equation}
\label{lemma_L-_above_L+_comparison_cond1_for_phi-}
\begin{aligned}
  G^{\mu}\left(v,\phi^-(v)\right)
  <
  \frac{\phi^-(v)}{v+\left(\phi^-(v)\right)^2\left(f(v)\right)^{-1}}
  &=
  \frac{\phi^-(v)f(v)}{vf(v)+\left(\phi^-(v)\right)^2},
  \\
  &\forall v\in\left(v_1^\mu,v_2^\mu\right)
  \cup\left(v_2^\mu,v_2^0\right)
  \cup\left(v_2^0,0\right).
\end{aligned}
\end{equation}
In terms of
\eqref{lemma_L-_above_L+_initial_cond1_for_phi-},
\eqref{lemma_L-_above_L+_comparison_cond1_for_phi-}
and applying Theorem \ref{thm_comparison}
and Remark \ref{remark_comparison}
on
$
\left(v_1^\mu,v_2^\mu\right),
$
we obtain
$\varphi^{\mu}(v)<\phi^-(v)$
on
$
\left(v_1^\mu,v_2^\mu\right).
$
Then, since
$\phi^-(v)$ and
$\varphi^{\mu}(v)$
continuously depend on $v$
and
$\phi^-(v)$
strictly increases for $v\in$
$
\left[v_2^\mu,v_2^0\right],
$
we have
\begin{equation}
\label{lemma_L-_above_L+_initial_cond2_for_phi-}
  \varphi^{\mu}(v_2^0)
  =
  \varphi^{\mu}(v_2^\mu)
  \le
  \phi^-(v_2^\mu)
  <
  \phi^-(v_2^0).
\end{equation}
In terms of
\eqref{lemma_L-_above_L+_comparison_cond1_for_phi-},
\eqref{lemma_L-_above_L+_initial_cond2_for_phi-}
and applying Theorem \ref{thm_comparison}
and Remark \ref{remark_comparison} on
$
\left(v_2^0,0\right),
$
we obtain
$\varphi^{\mu}(v)<\phi^-(v)$
for any
$
v\in\left(v_2^0,0\right),
$
and hence,
\begin{equation*}
  \varphi^{\mu}(0)\le\phi^-(0).
\end{equation*}
See Figure
\ref{fig:Possible cases for varphi and psi_varphi0=0}
and
\ref{fig:Possible cases for varphi and psi_varphi0>0}
for possible behavior of $\varphi^\mu(v).$

One the other hand, let
$\lambda\in(-1,0)$ and consider the following
piecewise differential equation
\begin{equation*}\label{piecewise diffeq+}
\frac{\mathrm{d}v}{\mathrm{d}z}
=H^{\lambda}(v,z)\triangleq
\left\{
       \begin{aligned}
       &\frac{v}{z}+\frac{z}{\lambda},  &z\in\left(0,z^{\lambda}_3\right],\\
       &\frac{v}{z}-z,          &z\in\left(z^{\lambda}_3,+\infty\right),
       \end{aligned}
\right.
\end{equation*}
where
\begin{equation*}
  z_3^{\lambda}=\sqrt{-\lambda}
  \left(
  \sqrt{v_3^{0}}+\sqrt{v_3^0-v_3^{\lambda}}
  \right),
\end{equation*}
with initial point $(v,z)=\left(v_3^{0},\sqrt{-\lambda v_3^{0}}\right).$
Its solution $v=\psi^{\lambda}(z)$
is explicitly formulated by
\begin{equation*}\label{explicit_sol_of_psilambda}
  \psi^{\lambda}(z)=\left\{
  \begin{aligned}
        &\frac{z^2}{\lambda}+2\sqrt{-\frac{v_3^{0}}{\lambda }}z,
        &z\in\left[\sqrt{-\lambda v_3^{0}},z^{\lambda}_3\right], \\
        &-z^2+\left(
        \frac{v_3^{\lambda}}{z_3^{\lambda}}
        +z_3^{\lambda} \right)z,
        &z\in\left(z^{\lambda}_3,+\infty\right).
  \end{aligned}
  \right.
\end{equation*}
Note that
$$
v_3^\lambda=\psi^\lambda(z_3^\lambda)
=\lim_{z\rightarrow z_3^\lambda,z>z_3^\lambda}
\psi^\lambda(z)
=\lim_{z\rightarrow z_3^\lambda}\left[
-z^2+\left(
        \frac{v_3^{\lambda}}{z_3^{\lambda}}
        +z_3^{\lambda} \right)z
\right],
$$
and that the equation
$
v_3^{\lambda}
=-z^2+\left(
        v_3^{\lambda}/z_3^{\lambda}
        +z_3^{\lambda} \right)z
$
has two real roots (counting with multiplicity) the
larger of which is denoted by $\bar z^{\lambda}$ and
satisfies
$
\sqrt{v_3^{\lambda}}
\le
\bar z^{\lambda}
<v_3^{\lambda}/z_3^{\lambda}+z_3^{\lambda}.
$
Letting
$
I^{\lambda}=
\left(\sqrt{-\lambda v_3^{0}},z^{\lambda}_3\right)
\cup\left(\bar z^{\lambda},
v_3^{\lambda}/z_3^{\lambda}+z_3^{\lambda}\right),
$
One can see that
$$
H^{\lambda}\left(\psi^\lambda(z),z\right)
=\frac{\mathrm{d}\psi^\lambda}{\mathrm{d}z}<0,
\quad\forall z\in I^\lambda.
$$
Letting $\omega^\lambda(v)$ be the inverse
of $\psi^{\lambda}(z)$ on $I^{\lambda},$
we have
\begin{equation}\label{lemma_L-_above_L+_initial_cond1_for_phi+}
\lim_{v\rightarrow v_3^0}\omega^\lambda(v)
=\sqrt{-\lambda v_3^0}
>0
=\lim_{v\rightarrow v_3^0}\phi^+(v)
\end{equation}
and
$
\lim_{v\rightarrow v_3^\lambda,v<v_3^\lambda}\omega^\lambda(v)
=\bar z^\lambda
\ge z_3^\lambda
=\lim_{v\rightarrow v_3^\lambda,v>v_3^\lambda}\omega^\lambda(v).
$
Noting that
\begin{equation*}
0>f(v)>\left\{
  \begin{aligned}
  &-1, &v\in\left(0,v^{\lambda}_3\right],\\
  &\lambda, &v\in\left(v^{\lambda}_3,v^{0}_3\right],
  \end{aligned}
\right.
\end{equation*}
we have for any
$
v\in \left(0,v_3^\lambda\right)
\cup\left(v_3^\lambda,v_3^0\right)
=\psi^\lambda\left(I^\lambda\right),
$
$$
\frac{v}{\omega^{\lambda}(v)}
+\frac{\omega^{\lambda}(v)}{f(v)}
=
\frac{\psi^\lambda(z)}{z}
+\frac{z}{f\left(\psi^\lambda(z)\right)}
<
H^{\lambda}\left(\psi^\lambda(z),z\right)
=
\frac{\mathrm{d}\psi^\lambda}{\mathrm{d}z}
<0,
$$
and hence,
\begin{equation}
\label{lemma_L-_above_L+_comparison_cond1_for_phi+}
\frac{\omega^{\lambda}(v)f(v)}
{vf(v,\alpha)+\left(\omega^{\lambda}(v)\right)^2}
>\frac{1}{H^{\lambda}\left(v,\omega^{\lambda}(v)\right)}
=\frac{\mathrm{d}\omega^{\lambda}}{\mathrm{d}v},
\quad
\forall v\in \left(0,v_3^\lambda\right)
\cup\left(v_3^\lambda,v_3^0\right).
\end{equation}
In terms of \eqref{lemma_L-_above_L+_initial_cond1_for_phi+},
\eqref{lemma_L-_above_L+_comparison_cond1_for_phi+}
and applying Theorem \ref{thm_comparison} and Remark \ref{remark_comparison} on
$\left(v_3^\lambda,v_3^0\right),$
we obtain
$
\phi^+(v)<\omega^\lambda(v)
$
on
$\left(v_3^\lambda,v_3^0\right),$
and hence,
\begin{equation}
\label{lemma_L-_above_L+_initial_cond2_for_phi+}
\phi^+(v_3^\lambda)
\le
\lim_{v\rightarrow v_3^\lambda,v>v_3^\lambda}\omega^\lambda(v)
=
z_3^\lambda
\le
\bar z^\lambda
=
\lim_{v\rightarrow v_3^\lambda,v<v_3^\lambda}\omega^\lambda(v).
\end{equation}
If
$
\phi^+(v_3^\lambda)=z_3^\lambda=\bar z^\lambda,
$
then
$$
\begin{aligned}
\frac{v_3^\lambda}{\phi^+(v_3^\lambda)}
+\frac{\phi^+(v_3^\lambda)}{f(v_3^\lambda)}
=
\frac{v_3^\lambda}{z_3^\lambda}
+\frac{z_3^\lambda}{f(v_3^\lambda)}
<
\frac{v_3^\lambda}{z_3^\lambda}-z_3^\lambda
=
\lim_{z\rightarrow z_3^\lambda,z>z_3^\lambda}
H\left(\psi^\lambda(z),z\right)
\le
0.
\end{aligned}
$$
Hence,
\begin{equation}
\label{lemma_L-_above_L+_comparison_cond2_for_phi+}
\begin{aligned}
\left.
\frac{\partial\phi^+}{\partial v}
\right|_{v=v_3^\lambda}
>&
\lim_{z\rightarrow z_3^\lambda,z>z_3^\lambda}
\frac{1}{H\left(\psi^\lambda(z),z\right)}
=
\lim_{v\rightarrow v_3^\lambda,v<v_3^\lambda}
\frac{1}{H^{\lambda}\left(v,\omega^{\lambda}(v)\right)}
\\
=&
\lim_{v\rightarrow v_3^\lambda,v<v_3^\lambda}
\frac{\mathrm{d}\omega^\lambda}{\mathrm{d}v}
\quad\text{if}\quad
\phi^+(v_3^\lambda)=z_3^\lambda=\bar z^\lambda.
\end{aligned}
\end{equation}
In terms of
\eqref{lemma_L-_above_L+_comparison_cond1_for_phi+},
\eqref{lemma_L-_above_L+_initial_cond2_for_phi+},
\eqref{lemma_L-_above_L+_comparison_cond2_for_phi+}
and applying Theorem \ref{thm_comparison} and Remark \ref{remark_comparison} on
$\left(0,v_3^\lambda\right),$
we obtain
$
\phi^+(v)<\omega^\lambda(v)
$
for any
$v\in\left(0,v_3^\lambda\right),$
and hence,
$$
\phi^+(0)
\le
\lim_{v\rightarrow 0}\omega^\lambda(v)
=z_3^{\lambda}+\frac{v_3^{\lambda}}{z_3^{\lambda}}.
$$
See Figure
\ref{fig:Possible cases for varphi and psi_z3=barz}
and
\ref{fig:Possible cases for varphi and psi_z3<barz}
for possible behavior of
$\psi^\lambda(z).$

\begin{figure}[htbp]
\centering
\subfigure[
$\varphi^\mu(0)=0\ge z_2^\mu+v_2^0/z_2^\mu$
]{
\label{fig:Possible cases for varphi and psi_varphi0=0}
\includegraphics[width=0.47\linewidth,height=3cm]
{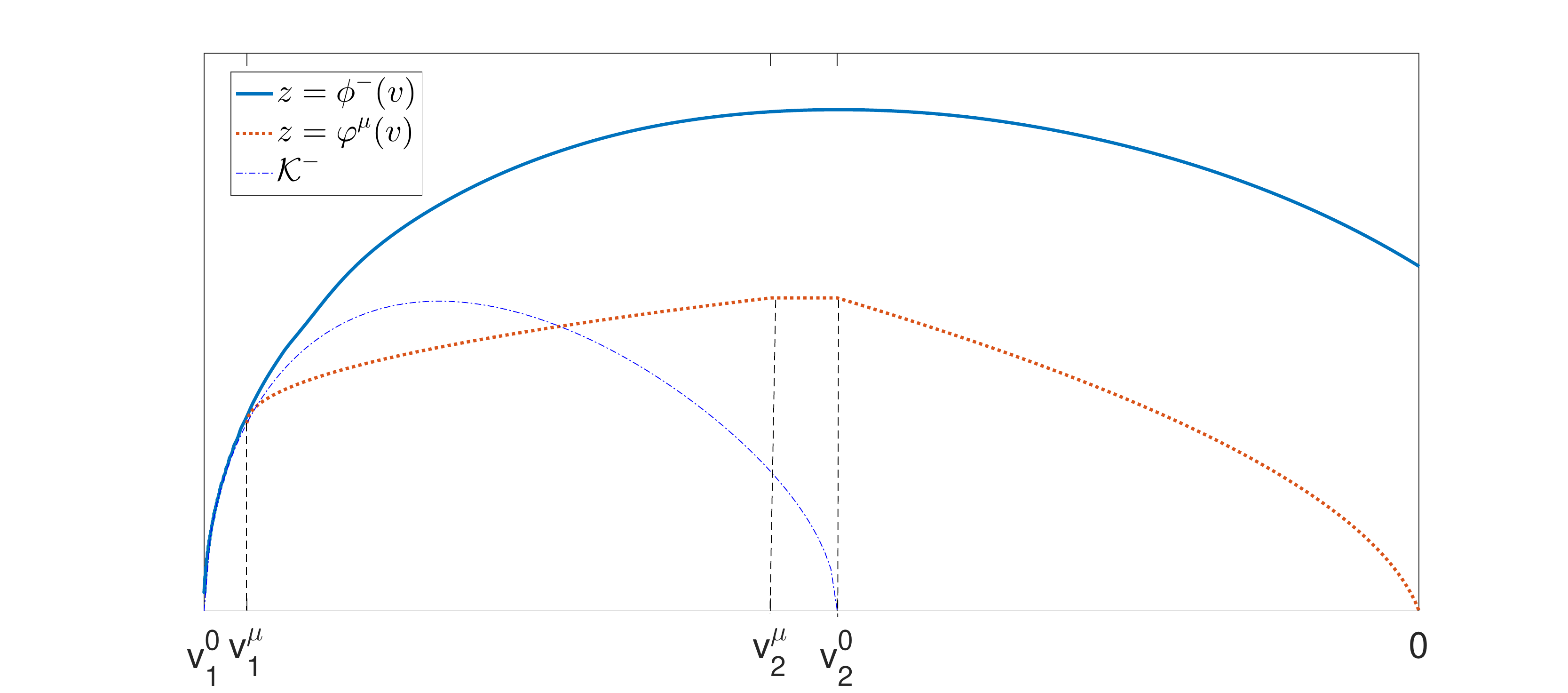}
}
\subfigure[
$\varphi^\mu(0)=z_2^\mu+v_2^0/z_2^\mu>0$
]{
\label{fig:Possible cases for varphi and psi_varphi0>0}
\includegraphics[width=0.47\linewidth,height=3cm]
{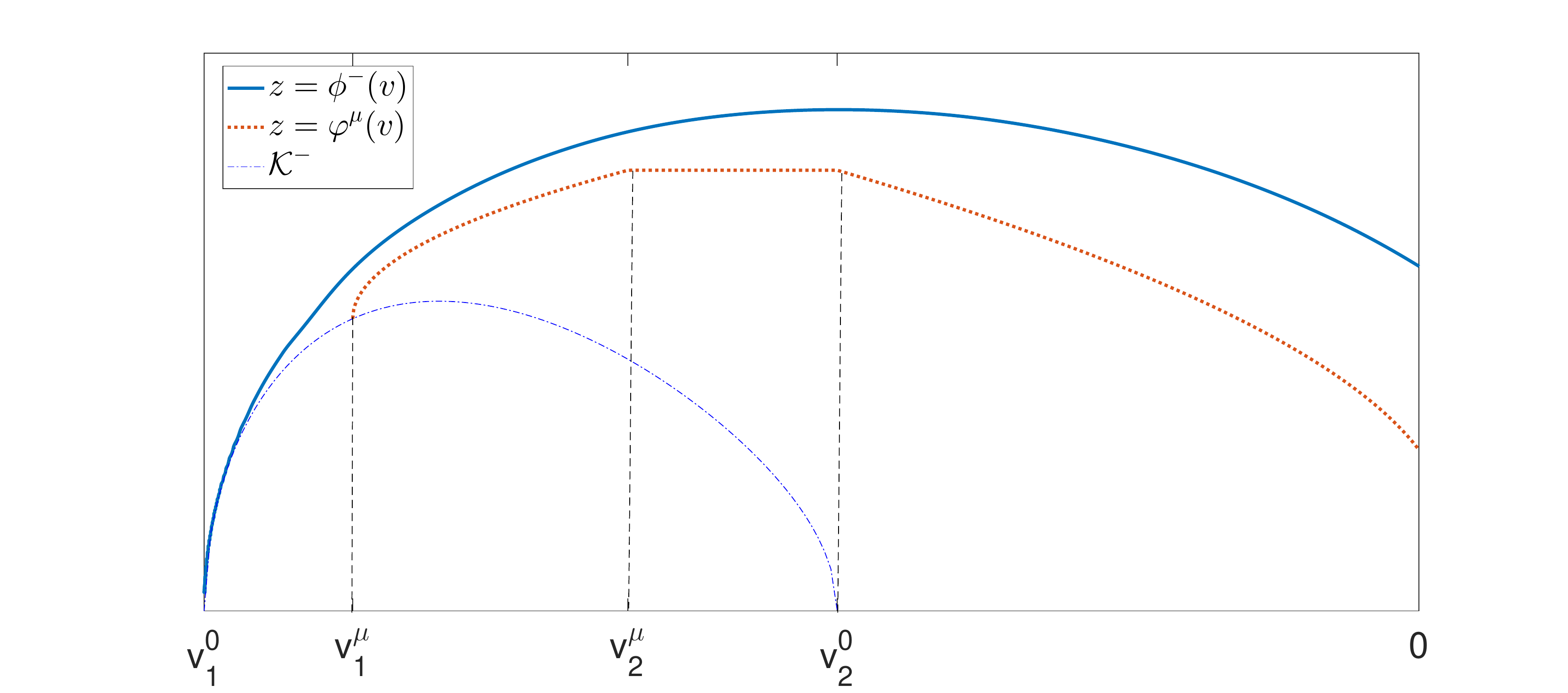}
}
\\
\subfigure[
$z_3^\lambda=\bar z$
]{
\label{fig:Possible cases for varphi and psi_z3=barz}
\includegraphics[width=0.47\linewidth,height=3cm]
{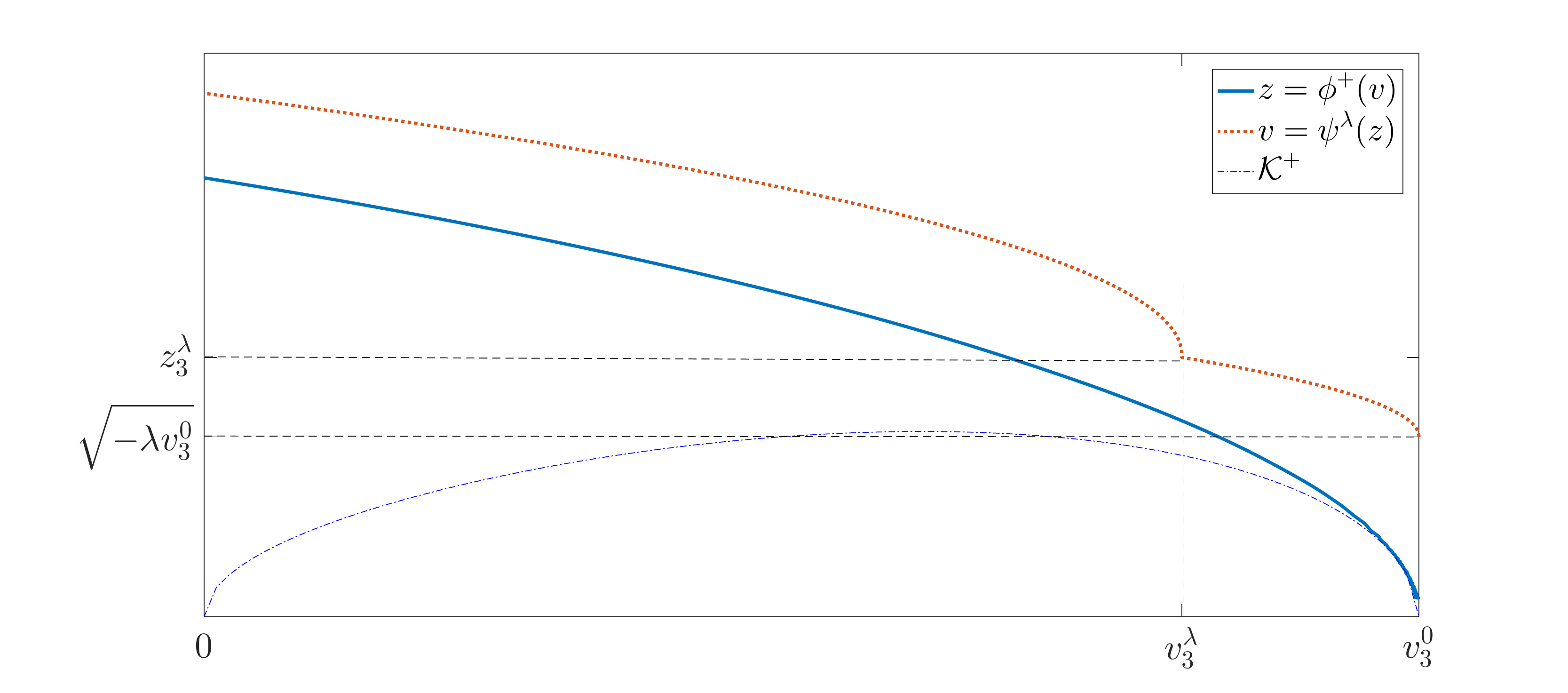}
}
\subfigure[
$z_3^\lambda<\bar z$
]{
\label{fig:Possible cases for varphi and psi_z3<barz}
\includegraphics[width=0.47\linewidth,height=3cm]
{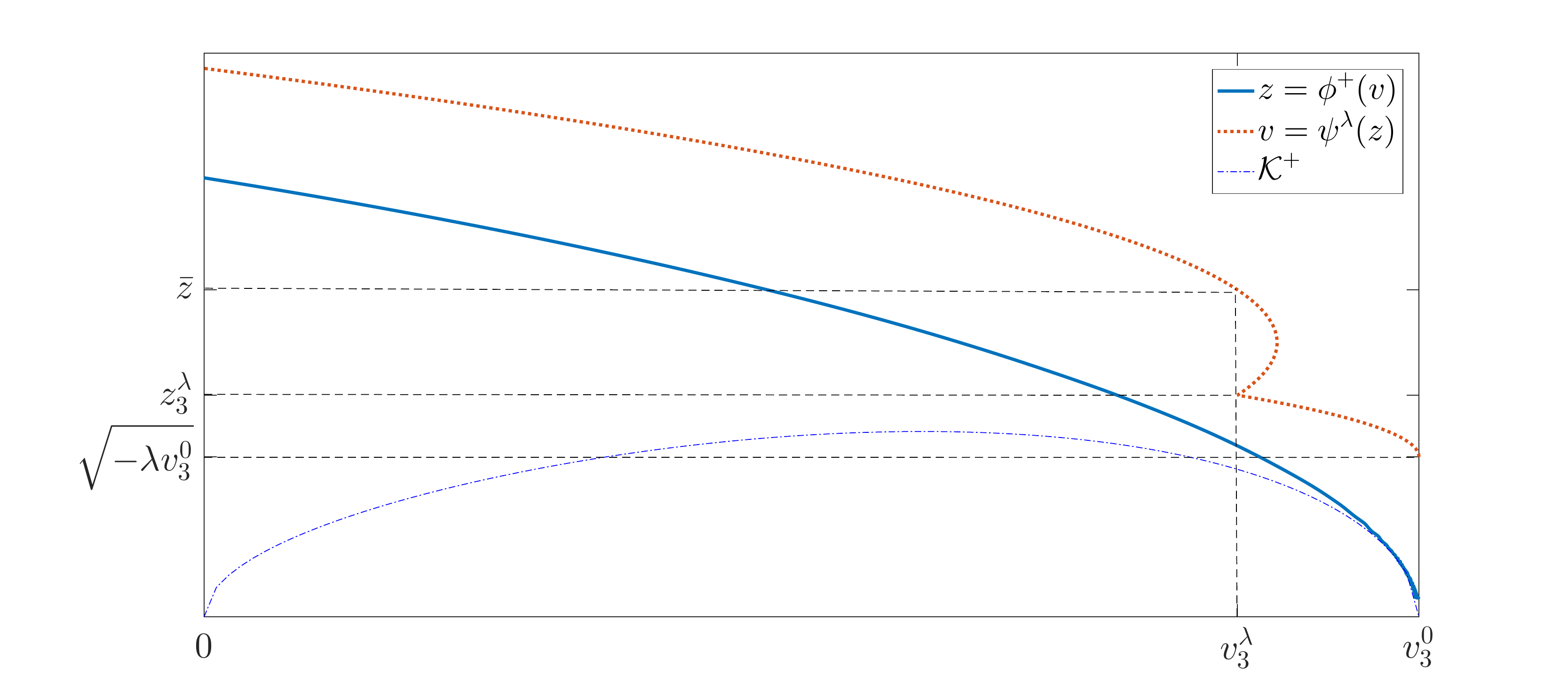}
}
\caption{Possible cases for
$\varphi^\mu$ and $\psi^\lambda$ with different
$\mu$ and $\lambda,$ respectively.}
\label{fig:Possible cases for varphi and psi}.
\end{figure}

Now from the structures of $\varphi^{\mu}(v)$ and $\psi^{\lambda}(z)$ ($\omega^\lambda(v)$),
we see that
\begin{equation}
\label{lemma_L-_above_L+_estimation_of_phi}
  \begin{aligned}
  \phi^-(0;\alpha)
  \ge&
  \varphi^{\mu}(0)
  =
  \max\left\{
  z_2^{\mu}+\frac{v_2^{0}}{z_2^{\mu}},0\right\} \\
  \ge&
  \max\left\{
  \sqrt{\mu}\left(\sqrt{-v_1^{\mu}}
  +\sqrt{v_2^{\mu}-v_1^{\mu}}\right)
  +\frac{v_2^{0}}{\sqrt{\mu}\left(\sqrt{-v_1^{\mu}}
  +\sqrt{v_2^{\mu}-v_1^{\mu}}\right)},0
  \right\},
  \\
    \phi^+(0;\alpha)
    \le&
    \omega^\lambda(0)
     =
     z_3^{\lambda} +\frac{v_3^{\lambda}}{z_3^{\lambda}}
     \\=&
     \sqrt{-\lambda}\left(\sqrt{v_3^{0}}+\sqrt{v_3^0-v_3^{\lambda}}\right)
     +\frac{v_3^{\lambda}}
     {\sqrt{-\lambda}\left(\sqrt{v_3^{0}}+\sqrt{v_3^0-v_3^{\lambda}}\right) }.
  \end{aligned}
\end{equation}
When
$\alpha=-3/\sqrt[3]{2},$
the roots of $f(v)=0$ are
$$
v_1^0\left(-\frac{3}{\sqrt[3]{2}}\right)
=-\frac{\sqrt{3}+1}{\sqrt[3]{2}},\quad
v_2^0\left(-\frac{3}{\sqrt[3]{2}}\right)
=-\frac{1}{\sqrt[3]{2}},\quad
v_3^0\left(-\frac{3}{\sqrt[3]{2}}\right)
=\frac{\sqrt{3}-1}{\sqrt[3]{2}}.
$$
By choosing
$\mu=11/16$
($<f(2\alpha/3)=1$)
and
$\lambda=-9/16,$
we have
$$
v_1^{11/16}\left(-\frac{3}{\sqrt[3]{2}}\right)
=-\frac{3\sqrt{5}+3}{4\sqrt[3]{2}},
\quad
v_2^{11/16}\left(-\frac{3}{\sqrt[3]{2}}\right)
=-\frac{3}{2{\sqrt[3]{2}}},
\quad
v_3^{-9/16}\left(-\frac{3}{\sqrt[3]{2}}\right)
=\frac{1}{2{\sqrt[3]{2}}}.
$$
Then from
\eqref{lemma_L-_above_L+_estimation_of_phi},
we have
\begin{equation*}
  \begin{aligned}
  \phi^-\left(0;-\frac{3}{\sqrt[3]{2}}\right)
  &\ge\left.
  \sqrt{\mu}\left(\sqrt{-v_1^{\mu}}
  +\sqrt{v_2^{\mu}-v_1^{\mu}}\right)
  +\frac{v_2^{0}}{\sqrt{\mu}\left(\sqrt{-v_1^{\mu}}
  +\sqrt{v_2^{\mu}-v_1^{\mu}}\right)}
  \right|_{\mu=11/16,\alpha=-3/\sqrt[3]{2}}
  \\
  &=\sqrt[3]{2}\sqrt{33}
  \left(\frac{\sqrt{\sqrt{5}+2}}{8}
  -\frac{4}{33}\sqrt{\sqrt{5}-2}\right)
  \approx1.4358
  \end{aligned}
\end{equation*}
and
\begin{equation*}
  \begin{aligned}
  \phi^+\left(0;-\frac{3}{\sqrt[3]{2}}\right)
  &\le
  \left.
  \sqrt{-\lambda}\left(\sqrt{v_3^{0}}+\sqrt{v_3^0-v_3^{\lambda}}\right)
     +\frac{v_3^{\lambda}}
     {\sqrt{-\lambda}\left(\sqrt{v_3^{0}}+\sqrt{v_3^0-v_3^{\lambda}}\right) }
  \right|_{
     \lambda=-9/16,\alpha=-3/\sqrt[3]{2}
  } \\
  &=\left(\frac{25}{12}-\frac{7}{24}\sqrt{3-\sqrt{3}}\right)
  \frac{\sqrt{\sqrt{3}-1}}{\sqrt[6]{2}}
  \approx1.3377.
  \end{aligned}
\end{equation*}
Hence,
\begin{equation}\label{lemma_L-_above_L+_alpha=2.3737}
\phi^-\left(0;-\frac{3}{\sqrt[3]{2}}\right)
>
\phi^+\left(0;-\frac{3}{\sqrt[3]{2}}\right).
\end{equation}

\end{proof}
\end{lem}

\section{Uniqueness and existence range of the limit cycle}

We first establish the relation between
the limit cycle of system \eqref{concretesys}
and the behavior of the separatrices of system \eqref{transformedSys(v,z)}.

\begin{figure}[htbp]
\centering
\subfigure[$\phi^-(0;\alpha)<\phi^+(0;\alpha)$]{
\label{fig:phase_portraits_on_disc_phi-(0)<phi+(0)}
\includegraphics[width=0.45\linewidth]{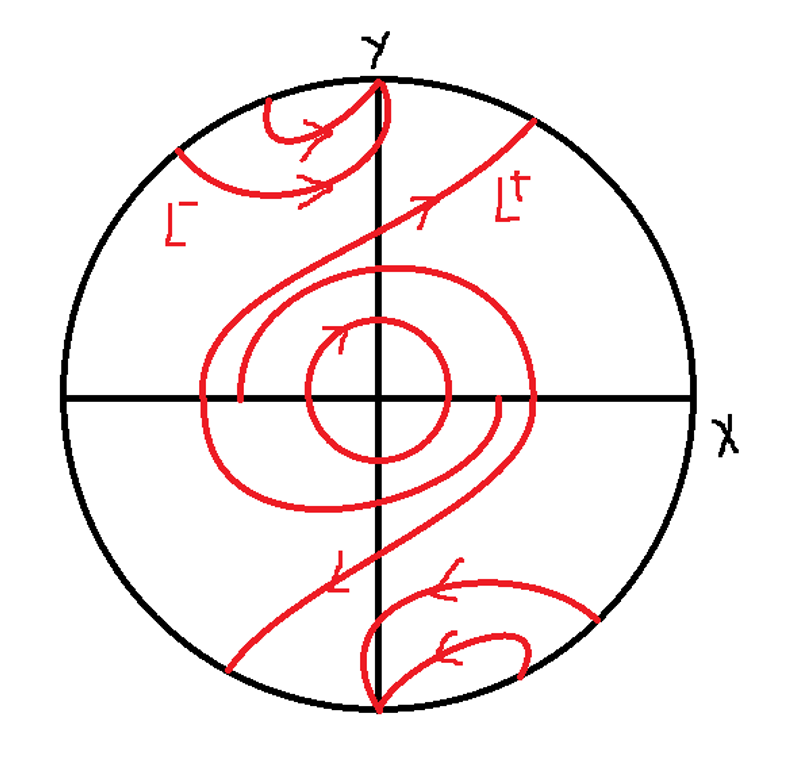}
}
\subfigure[$\phi^-(0;\alpha)>\phi^+(0;\alpha)$]{
\label{fig:phase_portraits_on_disc_phi-(0)>phi+(0)}
\includegraphics[width=0.45\linewidth]{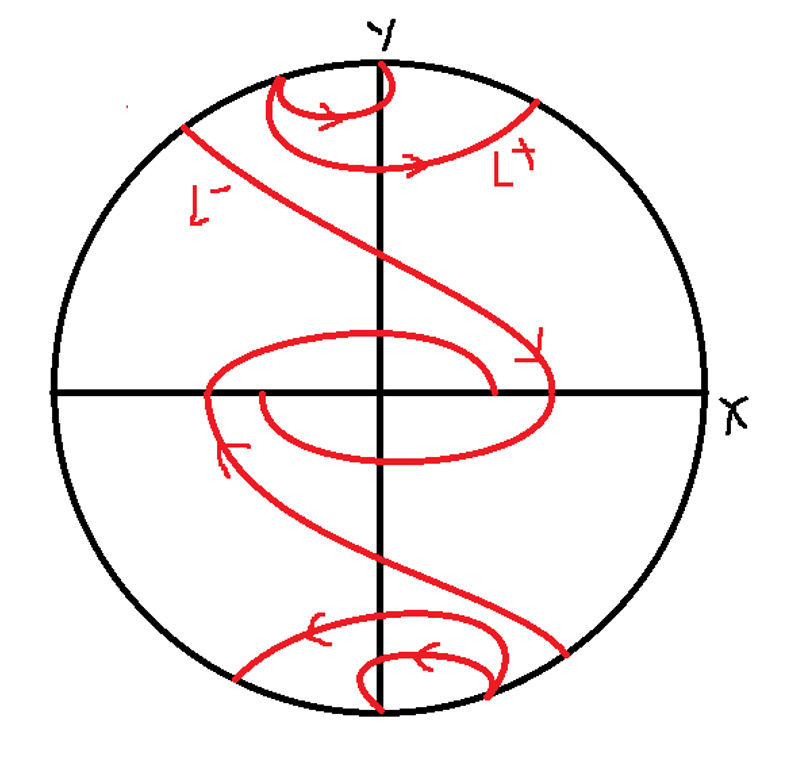}
}
\caption{Possible cases for the phase portraits of system
\eqref{concretesys} on the Poincar\'{e} disc}    
\label{fig:phase_portraits_on_disc}    
\end{figure}

\begin{lem}\label{lemma_relation_between_L+-_and_LC}
Assume that $\alpha<-3/\sqrt[3]{4}.$

(i) If the separatrices in the transferred
system \eqref{transformedSys(v,z)} satisfy
$\phi^-(0;\alpha)<\phi^+(0;\alpha),$
i.e., the case shown in
Figure \ref{fig:sys(v,z)_alpha=2}
or \ref{fig:sys(v,z)_alpha=2.15},
then system \eqref{concretesys} has at least one limit cycle.

(ii) If system \eqref{concretesys} has at least one limit cycle
and the outermost one is externally unstable,
then $\phi^-(0;\alpha)\le\phi^+(0;\alpha),$
i.e., the case shown in
Figure \ref{fig:sys(v,z)_alpha=2},
\ref{fig:sys(v,z)_alpha=2.15}
or \ref{fig:sys(v,z)_alpha=2.198}.

\begin{proof}
(i) From Proposition \ref{proposition_range_of_v+-}
we have
$v^+=-\infty.$
For system \eqref{transformedSys(v,z)}
in the $(v,z)$ plane,
$L^+$ positively approaches $P_3$
along the line $v=v_3^0,$
crosses positive $z$-axis and
negatively approaches infinity
in the second quadrant.
Then for system \eqref{concretesys}
in the $(x,y)$ plane,
$L^+$ positively approaches infinity
along the line $y=v_3^0 x$ in the upper half plane,
crosses positive $y$-axis, and negatively approaches
negative $x$-axis (see Figure
\ref{fig:phase_portraits_on_disc_phi-(0)<phi+(0)}).
Since $(\pm\infty,0)$ are not critical points of
system \eqref{concretesys}
and the system is central symmetry,
$L^+$ must spiral outward the origin.
On the other hand, the origin is a stable focus.
Therefore, the $\alpha$-limit set of $L^+$
must be a limit cycle.

(ii) Assume that
 $\phi^-(0;\alpha)>\phi^+(0;\alpha)$
though the system has at least one limit cycle
and the outermost one is externally unstable.
From Proposition \ref{proposition_range_of_v+-}
we have
$v^-=+\infty.$
For system \eqref{transformedSys(v,z)}
in the $(v,z)$ plane,
$L^-$ negatively approaches $P_1$
along the line $v=v_1^0,$
crosses positive $z$-axis and
positively approaches infinity
in the first quadrant.
Then for system \eqref{concretesys}
in the $(x,y)$ plane,
$L^-$ negatively approaches infinity
along the line $y=v_1^0 x$ in the upper half plane,
crosses positive $y$-axis, and positively approaches
positive $x$-axis (see Figure
\ref{fig:phase_portraits_on_disc_phi-(0)>phi+(0)}).
Since $(\pm\infty,0)$ are not critical points of
system \eqref{concretesys}
and the system is central symmetry,
$L^-$ must spiral toward the outermost
and externally unstable
limit cycle, which is a contradiction. Hence,
$\phi^-(0;\alpha)\le\phi^+(0;\alpha).$
\end{proof}
\end{lem}

Then with the help of known results,
we immediately obtained that system
\eqref{transformedSys(v,z)}
has a unique heteroclinic separatrix
in the upper half plane when $\alpha$ varies.
\begin{thm}\label{thm_unique_heteroclinic_separatrix}
There is a unique $\alpha^* \in \left(-\infty,-3/\sqrt[3]{4}\right)$
such that
$\phi^-(0;\alpha^*)=\phi^+(0;\alpha^*),$
i.e. the case shown in Figure
\ref{fig:sys(v,z)_alpha=2.198}.
Furthermore, we have
$$
\frac{-3}{\sqrt[3]{2}}\,(\approx-2.3811)
<\alpha^*\le
\frac{-3\sqrt[3]{6\sqrt3-9}}{\sqrt[3]4}
\,(\approx-2.1103).
$$
\begin{proof}
Recall that from Lemma \ref{lemma_Gasull2006}(ii)
and Lemma \ref{lemma_Giacomini2015},
system \eqref{concretesys}
has a unique and unstable limit cycle when
$0>\alpha>-3\sqrt[3]{6\sqrt3-9}/\sqrt[3]4
\approx -2.1103.$
Since Lemma \ref{lemma_alpha_continuity_of_phi-+}
gives the continuity of $\phi^-$ and $\phi^+$
with respect to $\alpha,$
then from Lemma \ref{lemma_relation_between_L+-_and_LC}(ii)
we have
$$
\phi^-\left(0;-3\sqrt[3]{\frac{6\sqrt3-9}{4}}\right)
\le\phi^+\left(0;-3\sqrt[3]{\frac{6\sqrt3-9}{4}}\right).
$$
Together with Lemma \ref{lemma_L-_above_L+}
which says
$$
\phi^-\left(0;-\frac{3}{\sqrt[3]{2}}\right)
>\phi^+\left(0;-\frac{3}{\sqrt[3]{2}}\right),
$$
there is a
$\alpha^*\in\left(
-3\sqrt[3]{2},
-3\sqrt[3]{6\sqrt3-9}/\sqrt[3]4
\right]
$
such that
$$
\phi^-\left(0;\alpha^*\right)
=\phi^+\left(0;\alpha^*\right).
$$
The uniqueness of $\alpha^*$
on $\left(-\infty,-3/\sqrt[3]{4}\right)$
is deduced from Lemma \ref{lemma_L_monotonically_varies}.
\end{proof}
\end{thm}

\begin{rmk}
\normalfont{
By a numerical computation, we have
$\alpha^*\approx-2.198.$
}
\end{rmk}

At last, we give the uniqueness of the limit cycle
of system \eqref{concretesys}
after introducing two lemmas.

\begin{lem}[{\cite[Lemma 4.1(iv)]{gasull2006upper}}]
\label{lemma_Gasull2006(2)}
If for some $\bar\alpha<0$ system \eqref{concretesys}
has no limit cycles, the same holds for any
$\alpha\le\bar\alpha.$
\end{lem}

Since system \eqref{concretesys} is a
semi-complete family of rotated vector fields
(mod $xy=0$) with respect to $\alpha,$
we have the following result which is deduced
from
\cite[Lemma 4 and Remark 1]{perko1993rotated}.
\begin{lem}
\label{lemma_Perko1993}
Given any negatively oriented, externally unstable limit
cycle $\Gamma(\alpha_0)$ of system \eqref{concretesys}
with $\alpha=\alpha_0,$
there exists an outer neighborhood $N$ of $\Gamma(\alpha_0)$
such that through each point of $N$ there
passes a limit cycle $\Gamma(\alpha)$ of
the system with $\alpha<\alpha_0.$
\end{lem}

\begin{thm}\label{thm_unique_LC}
System \eqref{concretesys} has a unique limit cycle when
$\alpha\in(\alpha^*,0),$
where
$\alpha^*$ is given in Theorem
\ref{thm_unique_heteroclinic_separatrix},
while it has no limit cycles when
$\alpha\in(-\infty,\alpha^*]\cup[0,+\infty).$

\begin{proof}
Since from Lemma \ref{lemma_Giacomini2015}
the system has a limit cycle for
$
\alpha\in\left(-3\sqrt[3]{6\sqrt3-9}/\sqrt[3]4,0\right),
$
we can denote $\bar\alpha$ by
\begin{equation*}\label{thm_unique_LC_critical alpha}
  \bar\alpha=
  \inf\left\{\alpha \;| \;
  \text{system \eqref{concretesys} has at least one limit cycle}
  \right\}.
\end{equation*}
Lemma \ref{lemma_Gasull2006(2)}
implies that the system has at least one limit cycle
for any $\alpha\in(\bar\alpha,0).$
Apparently, there are three cases:
$\bar\alpha>\alpha^*,\;\bar\alpha<\alpha^*$
and
$\bar\alpha=\alpha^*.$
Recall that from Theorem
\ref{thm_unique_heteroclinic_separatrix},
$\alpha^*$ is unique
and satifies
$
-3/\sqrt[3]{2}\,(\approx-2.3811)
<\alpha^*
\le-3\sqrt[3]{6\sqrt3-9}/\sqrt[3]4
\,(\approx-2.1103).
$

If $\bar\alpha>\alpha^*,$ take
$
\beta\in
\left(\alpha^*,\bar\alpha\right).
$
From Lemma \ref{lemma_L_monotonically_varies},
we have
$\phi^-(0;\beta)<\phi^-(0;\alpha^*)
=\phi^+(0;\alpha^*)<\phi^+(0;\beta).$
Then form Lemma \ref{lemma_relation_between_L+-_and_LC}(i),
system \eqref{concretesys} with $\alpha=\beta$
has at least one limit cycle,
which contradicts to the
definition of $\bar\alpha.$

If $\bar\alpha<\alpha^*,$ take
$
\beta\in
\left(\max\left\{\bar\alpha,-3/\sqrt[3]{2}\right\},
\alpha^*\right).
$
From Lemma \ref{lemma_L_monotonically_varies},
we have
$\phi^-(0;\beta)>\phi^-(0;\alpha^*)
=\phi^+(0;\alpha^*)>\phi^+(0;\beta).$
From Lemma \ref{lemma_Gasull2006}(ii)
and the definition of $\bar\alpha,$
system \eqref{concretesys} with $\alpha=\beta$
has a unique and unstable limit cycle. Then
from Lemma \ref{lemma_relation_between_L+-_and_LC},
we have
$\phi^-(0;\beta)\le\phi^+(0;\beta),$
which is a contradiction.

From the above discussions, we have
$\bar\alpha=\alpha^*.$
At last, we show that the system has no limit cycles
when $\alpha=\alpha^*.$
Otherwise, the limit cycle, denoted by $\Gamma,$
must be unique and unstable
from Lemma \ref{lemma_Gasull2006}(ii),
and it is negatively oriented.
From Lemma \ref{lemma_Perko1993},
system \eqref{concretesys} with
some $\alpha<\alpha^*$
has a limit cycle passing through
some outer neighborhood of $\Gamma,$
which contradicts to the definition of
$\bar\alpha=\alpha^*.$

Since from Lemma \ref{lemma_Gasull2006}(i)
system \eqref{concretesys} has
no limit cycles when $\alpha\ge0,$
it has a unique limit cycle when
$0>\alpha>\bar\alpha=\alpha^*\approx-2.198,$
while it has no limit cycles for the else region.

\end{proof}
\end{thm}

\biboptions{numbers,sort&compress}

\end{document}